\documentclass[12pt]{article}
\usepackage[a4paper,margin=25mm]{geometry}
\usepackage[T1]{fontenc}
\usepackage{lmodern,amsmath,amssymb,amsthm,booktabs,graphicx,microtype,hyperref}
\hypersetup{hidelinks,pdfauthor={Kaizhe Zhang},pdfsubject={School of Energy and Environment, City University of Hong Kong, Hong Kong SAR, China},pdftitle={How much rate heterogeneity can finite kinetic data hide?}}
\newtheorem{theorem}{Theorem}
\newtheorem{proposition}[theorem]{Proposition}
\newtheorem{corollary}[theorem]{Corollary}
\newtheorem{lemma}{Lemma}
\newtheorem{supplementtheorem}{Theorem}
\newcommand{\res}{\mathcal R}
\allowdisplaybreaks[2]
\title{\bfseries How much rate heterogeneity can finite kinetic data hide?}
\author{Kaizhe Zhang\\[3pt]\small School of Energy and Environment\\
\small City University of Hong Kong, Hong Kong SAR, China\\
\small ORCID: \href{https://orcid.org/0009-0007-5676-122X}{0009-0007-5676-122X}\\
\small\href{mailto:kaizzhang9-c@my.cityu.edu.hk}{kaizzhang9-c@my.cityu.edu.hk}}
\date{}
\begin{document}
\maketitle
\begin{abstract}
How much signal can an unobserved kinetic rate carry while a positive mixture fits finitely many decay measurements? For readings in a closed error box we characterize the largest compatible mass on a closed rate set and distinguish attainment at finite rates from a limit involving an infinitely fast component. When a query and $r$ positive reference rates approach one interior rate, we identify a participating mass $\omega$ and an error scale $\eta_*$. The maximal queried mass divided by $\omega$ is $\xi/(1+\xi)+O(H)$, where $\xi=\eta/\eta_*$ stays in a compact subset of $(0,\infty)$ and $H$ is the cluster radius. The estimate is uniform over relative spacings and positive weights, even as either vanishes. For $2r$ consecutive readings, a small-coupling condition yields the unique physical quadrature at the alternating corner of the raw error box and a signed finite-separation correction with a uniform quadratic remainder. In a measured two-dye decay, the higher-count histogram requires at least 54.49\% of the first retained gate's signal between 1 and 8 ns, while two lifetimes outside 2--6 ns reproduce the same bands.
\end{abstract}
\par\smallskip\noindent{\small The original text and figures are licensed under CC BY 4.0. Figure~\ref{fig:measured}, which includes adapted third-party data, is separately licensed under CC BY-SA 4.0.}\par\smallskip
\section{What can a decay experiment resolve?}
A sample containing two fluorescent dyes can give an excellent single-exponential decay curve. This happens in the mixed solution studied below. Separate measurements put the dye lifetimes at $0.48$ and $4.03$ ns \cite{kaye}, yet a lifetime of $3.85$ ns fits all nine ratios in the higher-count histogram. The same readings do establish where some of the light came from. Every compatible mixture assigns more than half of the first retained gate's signal to lifetimes between 1 and 8 ns. They even admit a mixture with no signal between 2 and 6 ns. The precision needed to locate emission depends on the width of the question one asks.

A normalized decay is an average of exponentials over an unknown positive rate distribution. A fit chooses one explanation of the data. The experiment restricts a larger set: every distribution whose predictions lie in all measured bands. Maximizing the mass of a specified rate or range over this set asks which components can be excluded at a given precision. Minimizing it asks whether a component in that range is required. For gated fluorescence the mass is a fraction of the retained signal; conversion to molecular abundances also requires brightness and detection calibration.

\paragraph{From a particular experiment to a resolution law.}
For measured data and a chosen window, start with Theorem~\ref{thm:face}. For a proposed design with specified reference components and one queried rate, start with the two transition laws. Theorem~\ref{thm:face} takes the actual readings, their error bands and a target range. It gives matching ways to establish an extremum: construct a compatible population and combine the readings to bound every population. It applies at finite precision. The next two theorems ask how precision must improve as rates approach one another. Theorem~\ref{cl:nested} concerns one additional rate near a cluster. It gives the participating signal fraction and the limiting error scale at which half of that fraction can move to the query. Once mass and error are expressed on those two scales, the same curve describes the transition for each fixed sampling design. Theorem~\ref{thm:multicenter} treats consecutive sampling times. Under its coupling condition, it identifies a specific set of extreme readings and the unique positive mixture that attains the maximum. Its signed correction quantifies the effect of rates lying on opposite sides of the query.

\paragraph{Relation to earlier work.}
Whittall and MacKay constructed extreme NMR models, Parker and Song bounded spectral functionals, and Watson studied optimal positive-exponential fits \cite{whittall,parker,watson}. Classical moment theory gives extremal interval probabilities and maximum point masses \cite{ks,bp,simon,hurlimann}. Kre\u{\i}n and Nudel'man also treated rectangular moment uncertainty \cite[Chapter VI]{krein}. Positivity and finite representation of truncated moments play a central role in these calculations \cite{cf,didio}.

For a prescribed point mass added to a known measure, a normalized kernel change already contains the factor $u/(1+u)$ \cite[Theorem 4.2, Eq.~(37) and Example 4.1]{kernelpert}. Here the perturbation is specified by error bars on the original exponential readings. Theorem~\ref{cl:nested} derives the precision scale from those bars and constructs a finite physical mixture that reaches the leading bound, uniformly over relative gaps and positive weights that may vanish in the limit. Theorem~\ref{thm:multicenter} identifies the maximizing corner of the raw reading box and gives its signed correction with a uniform quadratic remainder. Results on colliding Prony and Fourier sources quantify sensitivity of finite lists of locations and amplitudes \cite{prony,fourier,liu}. Denoyelle et al.\ study positive-measure support recovery \cite{denoyelle}. Eftekhari et al.\ control every approximately feasible nonnegative measure in generalized Wasserstein distance and bound mass near separated sites under certificate hypotheses \cite[Theorems 11--12]{eftekhari}. Our query is a discontinuous point-mass functional, for which a transport bound alone supplies no sharp mass extremum. Inverse Vandermonde estimates make spacing dependence explicit in polynomial interpolation \cite{vander}. Related statistical work studies fitted parameters in identifiable mixtures \cite{le}. The quantity optimized here is the mass at a specified rate across all compatible positive distributions.

A companion study by Zhang and Zhao \cite{companion} considers a fixed
queried mass and minimizes a covariance-weighted observation distance over
the physical rate axis. Here the error set is a deterministic box in the
original readings, and the mass is optimized. Their Theorem~6 inverts a distance law for both weighted boxes and covariance ellipsoids. On the common exponential box model it recovers the transition in our earlier preprint, with their separation constant $\gamma=1/C$ in our normalization. Their Lemmas A and B credit the earlier local constructions. A covariance metric changes the feasible geometry and the separation constant. The present proofs and the
certificates for finite reading bands are self-contained.

Section~\ref{sec:finite} shows how to read an upper and a lower certificate on a two-rate example. Section~\ref{sec:transition} derives the scales behind the resolution law. The controlled calculations test where the law is accurate at finite separation. Section~\ref{sec:illustrations} returns to measured traces and the choice of an extra observation.

\section{Reading the range of compatible fractions}\label{sec:finite}
Suppose an experimentalist wants to rule out a large, rapidly decaying component. Three inputs suffice: the sampling times, an interval of plausible readings at each time, and the range of rates that would count as fast. The output is the largest fraction of signal that could come from that range. The guaranteed lower bound in an open range is one minus the supremum in its closed complement. This is an infimum over physical distributions, which a finite-rate distribution can sometimes only approach. Endpoints belong to the complement in that calculation.

Let $\mu$ be a probability distribution on the physical rates $D=[0,\infty)$.
With time and rate in reciprocal units, a normalized remaining signal has the form
\begin{equation}\label{eq:observation}
 s_\mu(t)=\int_D e^{-tk}\,d\mu(k).
\end{equation}
Measure \eqref{eq:observation} at distinct positive times $t_1<\dots<t_n$, and record $L_i\le s_\mu(t_i)\le U_i$.
These bounds may be asymmetric, and equality permits an exact reading.
Intersect the bands at repeated times. A reading at zero tests $s_\mu(0)=1$.
Put $F(k)=(e^{-t_i k})_{i=1}^n$ and $B=\prod_i[L_i,U_i]$.
For a closed target $E\subset D$ define the \emph{resolution function}
\begin{equation}\label{eq:resolution}
 \res_B(E)=\sup\big\{\mu(E):\mu\text{ a probability on }D,\ \int F\,d\mu\in B\big\}.
\end{equation}
Feasibility is tested before interpreting a fraction. The empty target has value zero whenever feasible.

\paragraph{Can three perfect readings hide a large component?}
To see how the ambiguity arises, measure the fraction that survives a time $\tau$ and use $X=e^{-\tau K}$, so that $s_\mu(j\tau)=\mathbb E X^j$. Suppose the physical rate range is known to give $1/4\le X\le3/4$. Equal reference weights at $X=3/8,5/8$ yield $\mathbb E X=1/2$ and $\mathbb E(X-1/2)^2=1/64$. If the first reading is exact and the second has half-width $\eta$, any candidate putting mass $p$ at $X=1/2$ satisfies
\[
 (1/64-\eta)_+\le\mathbb E(X-1/2)^2\le(1-p)/16.
\]
This calculation uses the first two readings. Put the largest permitted $p$ at the midpoint and split the rest equally between the endpoints. This distribution attains both inequalities at the optimum. Thus $\res_B(\{\log(2)/\tau\})=1-16(1/64-\eta)_+$. At $\eta=0$ the answer is exactly $3/4$. If the tolerance is zero, both distributions are symmetric about $1/2$. Their shared first two moments then also determine their third, so a third exact reading at $3\tau$ leaves them indistinguishable. Their readings at $4\tau$ differ by $3/4096$. If the first four readings all equal the reference values exactly, then $\mathbb E[(X-3/8)^2(X-5/8)^2]=0$, and every possible population is supported on the two reference rates. An exact fourth reading thus resolves the ambiguity in this example.

\begin{figure}[htbp]\centering
\includegraphics[width=\linewidth]{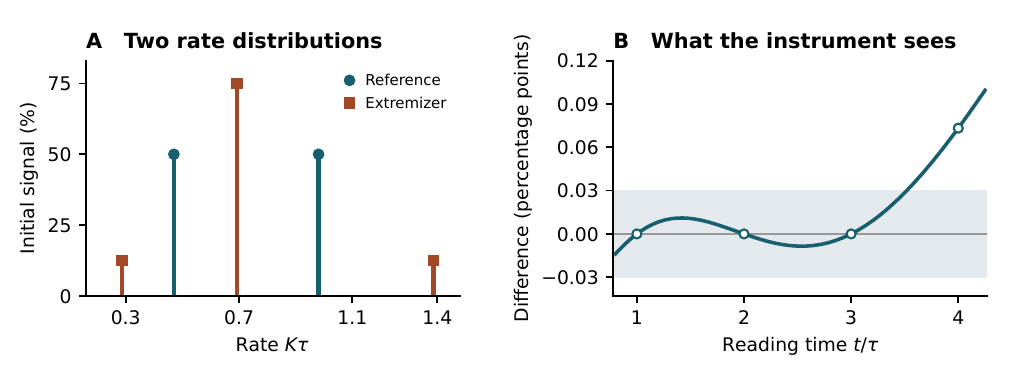}
\caption{Two populations with the same first three kinetic readings. The queried rate is $K\tau=\log 2$. Left, the bars show the exact amplitude fractions in the reference and in a population maximizing the queried fraction. Right, the signed difference of their normalized readings is magnified and expressed in percentage points. The shaded band is an illustrative $\pm0.03$ percentage-point allowance about the reference. The difference at the fourth reading exceeds twice this half-width, so the two possible observation intervals there are disjoint.}\label{fig:concept}
\end{figure}

\paragraph{How can a bound cover every possible rate?}
An upper certificate is a weighted sum of the measured channels. If its curve stays above one on the target and above zero everywhere else, averaging this inequality bounds the target mass of every compatible population. A second calculation exhibits rates and nonnegative weights that reach the bound or approach it with finite rates.

To include arbitrarily fast decays, adjoin the endpoint $k=\infty$ with $F(\infty)=0$. It represents an arbitrarily fast process with no signal at any positive measurement time. Write $\bar D=[0,\infty]$, $T=\overline E^{\,\bar D}$, $M=\operatorname{conv}F(\bar D)$ and $\bar{\res}$ for the compactified maximum. The convex hull $M$ is the set of physically limiting reading vectors. For $b\in\mathbb R^n$ set $h_B(b)=\sum_i(b_i^+U_i-b_i^-L_i)$ where $b_i^\pm=\max(\pm b_i,0)$.

\begin{theorem}[Finite-error resolution]\label{thm:face}
If $M\cap B\ne\varnothing$, the compactified value has two attained expressions,
\begin{equation}\label{eq:univ}
\begin{aligned}
 \bar{\res}
 &=\max_{\substack{k_j\in\bar D,\ w_j\ge0,\ \sum_jw_j=1\\\sum_jw_jF(k_j)\in B}}
       \sum_{j=1}^{n+1}w_j{\bf1}_T(k_j)\\
 &=\min_{\substack{c\in\mathbb R,\ b\in\mathbb R^n\\
          c+b^\top F(k)\ge{\bf1}_T(k)\ \forall k\in\bar D}}
        \{c+h_B(b)\}.
\end{aligned}
\end{equation}
If $B$ meets the interior of $M$, $\res_B(E)=\bar{\res}$, including for exact readings. In this case finite physical mixtures approach the value. It is attained physically precisely when a compact optimizer places zero mass at infinity. If $B$ meets only the boundary of $M$, it admits a unique compact law, whose mass at infinity determines physical feasibility. If $B\cap M$ is empty, no law is compatible.
\end{theorem}

Read the two lines of \eqref{eq:univ} in opposite directions. The first gives an example. Its weights are the fractions one could actually assign to the target, apart from the expressly identified limiting rate at infinity. The second gives a guarantee. A positive coefficient uses the upper end of its reading interval and a negative coefficient uses the lower end. When the example and guarantee agree, further fitting cannot change the bound. Rates with positive weight lie where the bounding curve touches the target indicator.

The proof in Appendix~\ref{app:finite} separates feasible reading vectors from larger target fractions. On the boundary, a nonnegative combination of exponential channels vanishes at all possible support points. Counting its zeros determines whether the optimizer uses the limiting fast rate. For instance, if the only reading is $\mathbb E e^{-K}=m\in(0,1)$, at most $m$ can lie at $K=0$. A mixture at $0$ and $R$ places $(m-e^{-R})/(1-e^{-R})$ there and approaches $m$ as $R$ grows. If a second exact reading gives $\mathbb E e^{-2K}=m$ as well, then $e^{-K}(1-e^{-K})$ must vanish almost surely. The two required sites are $K=0$ and $K=\infty$. This second experiment has no distribution on the physical finite-rate axis.

\paragraph{Using the result.}
For arbitrary sampling times, the numerical procedure proposes a discrete distribution and a bounding curve. A separate calculation checks the distribution against the original bands and the curve between grid points. When the box meets the interior of the moment body, a finite-rate repair makes the bracket converge under refinement (Supplement, Section~\ref{sec:grid}). Exact boundary readings can force a site absent from every mesh. The boundary construction in Appendix~\ref{app:finite} then identifies the support. For example, three readings at times $1/2,\sqrt{2},3$ give the finite-data brackets in Supplement, Section~\ref{app:practical}. Consecutive readings also admit a matrix calculation. For $b$ disjoint closed target intervals the sharp worst-case support count is $\min(n+1,\lfloor n/2\rfloor+2b+1)$ on a compact axis. Supplement, Section~\ref{irregular:section} proves this count and Supplement, Section~\ref{matrix:hausdorff} treats missing or singular moments.

The finite-error calculation can be repeated on a new dataset. What it does not by itself reveal is how the required precision changes when two reference rates approach each other. That is the subject of the next section.

\section{When does an extra rate become distinguishable?}\label{sec:transition}
Imagine two calibrated components whose rates are close. Could a third rate between them carry an appreciable fraction without changing any reading beyond its error bar? The answer depends on how much signal each reference component contributes, how far the rates are apart, and when and how accurately the instrument samples. The target in this section is a \emph{single specified rate}. A guarantee for a whole window uses the finite-data calculation in Section~\ref{sec:finite}.

Consider a reference population with $r\ge1$ distinct rates and positive weights $w_i$ summing to one. We query one additional rate, called the \emph{anchor}. Each measured channel may depart from the reference by $\eta\sigma_j$, where $\eta$ is an absolute signal tolerance and $\sigma_j$ describes relative precision. Denote the largest possible anchor fraction by $A$. All competing populations, including continuous distributions, enter its definition.

\paragraph{Three quantities to read first.}
The mass scale $\omega$ describes how much of the reference signal participates in the transition. The critical error $\eta_*$ is the precision at which roughly half of that mass can hide. The ratio $\xi=\eta/\eta_*$ places the actual experiment on this scale. At $\xi=1/3,1,3$ the leading prediction for the hidden share of $\omega$ is $1/4,1/2,3/4$. For example, if $\omega=0.20$ and $\xi=1$, the predicted hidden mass is about $0.10$ of the entire signal. A larger permitted error gives a larger compatible fraction. A weak component can make $\omega$ much smaller than one.

\paragraph{Computing the scales from the reference.}
In either of the coordinates specified below write the reference locations as $x_i$ and the anchor as $z\notin\{x_i\}$. Define
\begin{equation}\label{eq:cardinal}
 P(x)=\prod_{i=1}^r(x-x_i),\qquad D=P(z),\qquad
 \ell_i(x)=\frac{P(x)}{(x-x_i)P'(x_i)}
\end{equation}
with the quotient interpreted as a polynomial. The number $\ell_i(z)$ records how strongly reference site $i$ contributes when interpolating a polynomial of degree less than $r$ at the query. Thus the effect of a weak component depends on its position. The signed gap product $D=P(z)$ is zero precisely at a reference point. Here the query is distinct. Using \eqref{eq:cardinal}, define the mass scale $\omega$ and the signed normalization $d$ by
\begin{equation}\label{eq:scales}
 S=\sum_i\frac{\ell_i(z)^2}{w_i},\qquad
 \omega=S^{-1},\qquad d=D\sqrt{\omega},\qquad
 d^2=\left[\sum_i\frac{1}{w_i(z-x_i)^2 P'(x_i)^2}\right]^{-1}.
\end{equation}
Since $\sum_i\ell_i(z)=1$, one has $0<\omega\le1$. A small reference weight matters most when its interpolation coefficient is appreciable. The squared product in the last expression depends on the distance from each rate to the anchor and to every other reference rate.

For each sampling geometry below let $C$ be its specified design coefficient. Define
\begin{equation}\label{eq:master}
 \eta_*=\frac{d^2}{C},\qquad \xi=\frac{\eta}{\eta_*},\qquad
 f(\xi)=\frac{\xi}{1+\xi}.
\end{equation}
Thus $\xi$ compares experimental tolerance with the tolerance at which half the available mass scale $\omega$ can be hidden to leading order. The formulas below concern the transition where this ratio remains bounded away from zero and infinity. Other transitions can involve additional mass at larger tolerances.

\subsection{One cluster, with arbitrary fixed sampling times}\label{cl:section}
To predict the hardest part of a measurement, bring all the reference rates towards a common interior rate. Keep the instrument's sampling times fixed. Some rates may approach one another much faster than others, and a weak component may fade at the same time.

Write the physical reference rates as $k_*+x_i$ and the query as $k_*+z$ with fixed $k_*>0$. The cluster size $H=\max(|z|,\max_i|x_i|)$ tends to zero. Each configuration has distinct reference rates, and their gaps and positive weights may shrink at unrelated speeds. Fix $n\ge2r$ distinct positive times and fixed $\sigma_j>0$. Set $M_0=2r$, $t_0=0$, $g_j=e^{-k_*t_j}$ and
\begin{equation}\label{eq:design}
 C=\min_{|I|=M_0} M_0!\sum_{j\in I}
 \frac{\sigma_j}{g_j\left|\prod_{i\in I\cup\{0\},\,i\ne j}(t_j-t_i)\right|},
 \qquad I\subset\{1,\ldots,n\}.
\end{equation}
The minimization selects the most informative set of $2r$ readings for the leading bound. A change from rate coordinates $x$ to $a x$, with times divided by $a$, multiplies both $C$ and $d^2$ by $a^{2r}$. The normalized tolerance $\xi$ is dimensionless.

\paragraph{The scales in one calculation.}
Take one reference at $k_*=\log 2$ with weight one. Query $k_*+h$, initially at $h=0.1$, using times 1 and 2 with $\sigma_1=\sigma_2=1$. The reference displacement is zero, so $\ell_1=1$, $\omega=1$ and $d^2=h^2=0.01$. The reference readings are $g_1=1/2$ and $g_2=1/4$. The two time products in \eqref{eq:design} have absolute values 1 and 2. Hence
\[
 C=2!\left(\frac{1}{(1/2)\,1}+\frac{1}{(1/4)\,2}\right)=8,
 \qquad \eta_*=\frac{0.01}{8}=0.00125 .
\]
At $\eta=\eta_*$, the normalized error is $\xi=1$. If the displacement shrinks with $\eta=h^2/8$, the largest hidden fraction tends to one half. At the finite displacement $0.1$, the two-reading calculation gives $0.53787$ to five decimal places. The difference measures the finite-separation effect in this example.

\begin{theorem}[Resolution at the finest cluster scale]\label{cl:nested}
For fixed $r,k_*$ and design as above, uniformly over all configurations with $\xi$ in a fixed compact subinterval of $(0,\infty)$,
\begin{equation}\label{cl:formula}
                 A=\omega\big[f(\xi)+O(H)\big]\qquad (H\to0).
\end{equation}
The constant is independent of every relative gap and weight, in the fixed rate units. At most $r+1$ finite physical atoms, including the anchor, suffice for the lower estimate. The upper estimate also holds on the compactified rate axis.
\end{theorem}

Equation~\eqref{cl:formula} predicts a fraction of signal, in units of $\omega$. Dividing the error by $\eta_*$ puts different reference configurations on the same horizontal scale. Shrinking the entire cluster improves the approximation even when an internal gap or weight shrinks especially fast. For a fixed choice of sampling times and channel precisions, the normalized error bound works uniformly over gaps and weights on the stated range of $\xi$. Changing the design changes the constant.

One can read a critical exponent directly from the geometry. Suppose $w_i\asymp h^{\beta_i}$, $|z-x_i|\asymp h^{a_i}$ and $|x_i-x_j|\asymp h^{b_{ij}}$. Then
\[
 d^2\asymp h^\kappa,\quad
 \kappa=\max_i\big(\beta_i+2a_i+2\sum_{j\ne i}b_{ij}\big),\qquad
 \omega\asymp h^{\kappa-2\sum_i a_i}.
\]
Here $u\asymp v$ means their ratio stays bounded above and below by positive constants. Crowding can reduce the error needed to resolve a rate by many powers of its separation. If $\omega$ also shrinks, only a small fraction participates, even at the critical tolerance. Formula~\eqref{eq:scales} retains the constants that power counting alone loses.

With exactly $2r$ channels the maximizing distribution is unique for sufficiently small $H$. In increasing time order its readings reach the lower, upper, lower, upper ends of the bands in succession. With extra channels this exact conclusion holds on a minimizing subset if its alternating linearized perturbation lies strictly within the unused bands (Corollary~\ref{cl:corner}). This strict-slack condition is needed only for that exact subset conclusion. Theorem~\ref{cl:nested} uses the full design in all cases.

\paragraph{Why a rational transition appears.}
Take a polynomial with value one at the query. Its square bounds the mass at that point after averaging. Minimizing the mean square over polynomials of the allowed degree gives the Christoffel value. The minimizer vanishes at the other nodes of an attaining quadrature when the matching construction applies. In a basis scaled by the reference weights and the gap product, the leading quadratic form has an identity block for the reference population and a final entry $\xi$ supplied by the error allowance. Its minimum under the normalized evaluation constraint is $f(\xi)$. Constructing a positive physical distribution attaining the lower estimate and transferring the polynomial upper bound to the original exponential data are the substantive steps. Appendix~\ref{app:cluster} gives both.

\subsection{Separated clusters on a consecutive grid}\label{sec:multicenter}
Some experiments follow two groups that stay well separated, even as the rates inside each group come together. For equally spaced readings one can ask a stronger question. Which exact readings inside the error bars hide the most signal, and which mixture produces them?

Now assume the sampling times are $\tau,2\tau,\ldots,2r\tau$. Use the dimensionless survival coordinate $X=e^{-\tau K}$: the readings are moments $\mathbb E X^j$. In \eqref{eq:scales} the $x_i$ and $z$ now denote survival coordinates, all in a fixed interval $[a,b]\subset(0,1)$. Feasible laws may use all of $[0,1]$. Define $R_i(X)=P(X)/(X-x_i)$ and, with $[X^j]$ denoting the coefficient of $X^j$, put
\begin{equation}\label{eq:consecutive-design}
 C=\sum_{j=1}^{2r}\sigma_j|[X^j]P^2|,\quad
 C_i=\sum_{j=1}^{2r}\sigma_j|[X^j]R_iP|,\quad
 \pi_i=\frac{d^2}{w_i(z-x_i)^2 P'(x_i)^2}.
\end{equation}
Equations~\eqref{eq:scales} and \eqref{eq:consecutive-design} give $\sum_i\pi_i=1$. The coupling parameter
\begin{equation}\label{eq:coupling}
 \epsilon=\max_i |z-x_i|\sqrt{\pi_i}
\end{equation}
can be small even when several reference clusters stay a fixed distance apart. Both $\epsilon$ and $\xi$ can be computed without assigning the rates to clusters.

A quadrature here means a finite list of sites and positive weights that reproduces the measured moments. The theorem identifies exactly one such maximizing list under its hypotheses.

\begin{theorem}[Separated rates and the first correction]\label{thm:multicenter}
Fix $r,[a,b]$ as above, $0\le\sigma_j\le M_1$ and $\sigma_{2r}\ge c>0$. As $\epsilon\to0$, let $\xi$ stay in a fixed compact subinterval of $(0,\infty)$. For all sufficiently small $\epsilon$ the exact optimizer is the unique positive quadrature on the anchor and $r$ distinct additional points in $(0,1)$. Its raw moments are
\begin{equation}\label{eq:alternating-vertex}
 m_0^*=1,\qquad m_j^*=\sum_iw_i x_i^j+(-1)^j\eta\sigma_j,\quad 1\le j\le 2r .
\end{equation}
Writing $G_*=(m^*_{i+j})_{0\le i,j\le r}$ and $v(z)=(1,z,\ldots,z^r)^\top$, the exact mass is
\begin{equation}\label{eq:inverse-gram}
 A=\bigl(v(z)^\top G_*^{-1}v(z)\bigr)^{-1}.
\end{equation} With \eqref{eq:consecutive-design}--\eqref{eq:coupling}, its expansion is
\begin{equation}\label{eq:correction}
 \frac A\omega=f(\xi)-2 f(\xi)^2
       \sum_i \pi_i(z-x_i)\frac{C_i}{C}+O(\epsilon^2).
\end{equation}
The remainder constant depends only on $r,a,b,M_1,c$ and the fixed bounds for $\xi$.
\end{theorem}

Equation~\eqref{eq:alternating-vertex} selects the moment matrix in \eqref{eq:inverse-gram}. The optimizing reading is the lower band edge at the first time, the upper edge at the second, and so on. The matrix $G_*$ collects precisely these readings along its diagonals of constant index sum. Appendix~\ref{app:cluster} constructs its sites and proves that its bounding polynomial dominates every compatible distribution.

The signed displacement in \eqref{eq:correction} lets rates on opposite sides of the query contribute corrections of opposite signs. The total linear term is $O(\epsilon)$ since $\sum_i\sqrt{\pi_i}\le\sqrt r$ and $C_i/C$ is uniformly bounded. The zeroth-order law is therefore the same $f(\xi)$ as in Theorem~\ref{cl:nested}. On the common domain of the two theorems the change of coordinate gives $\xi_X/\xi_K=1+O(H)$ and $\omega_X/\omega_K=1+O(H)$; Proposition~\ref{lem:coordinates} proves these estimates uniformly even for vanishing weights. Inverting a raw numerical moment matrix near a collision can be unstable; the normalized form used in the proof gives a well-conditioned alternative when the individual reference parameters are available.

\paragraph{How large is the first correction?}
For two separated pairs with within-pair survival spacing $1/128$, the certified normalized maximum is $0.685900931$ to the displayed precision. The leading prediction is $0.686804748$. Including the signed correction gives $0.685905249$. The supplementary material specifies every rate, weight and tolerance. Section~\ref{sec:controlled8} compares exact finite-data answers with both predictions as the spacing changes.

\subsection{When does the common curve predict a finite experiment?}\label{sec:controlled8}
A measured sample usually has an unknown rate distribution. To test a formula for the largest hidden fraction, one needs a known reference and an independent calculation of that maximum. The following controls supply both. The reference rates and weights are prescribed. The reference centers are noiseless. In the instrument-based controls, the error widths come from the measured fluorescence or NMR records.

\paragraph{A common curve for nearby rates.}
For two equally weighted rates, start with displacements $0$ and $2\log(1+h)$ in fixed reciprocal-time units. Query the midpoint $\log(1+h)$, which carries zero reference mass. Reducing $h$ brings the rates together. Other controls reduce one weight to $h$, double the sampling interval, add a seventh-time reading, use a sparse pair of times, or adopt the relative widths of the fluorescence gates. The sparse pair uses one reference rate. For each design we choose five tolerances spanning $\xi=1/4$ to $4$. Theorem~\ref{cl:nested} predicts that dividing the allowed mass by $\omega$ brings these different experiments onto the same curve.

\begin{figure}[htbp]\centering
\includegraphics[width=\linewidth]{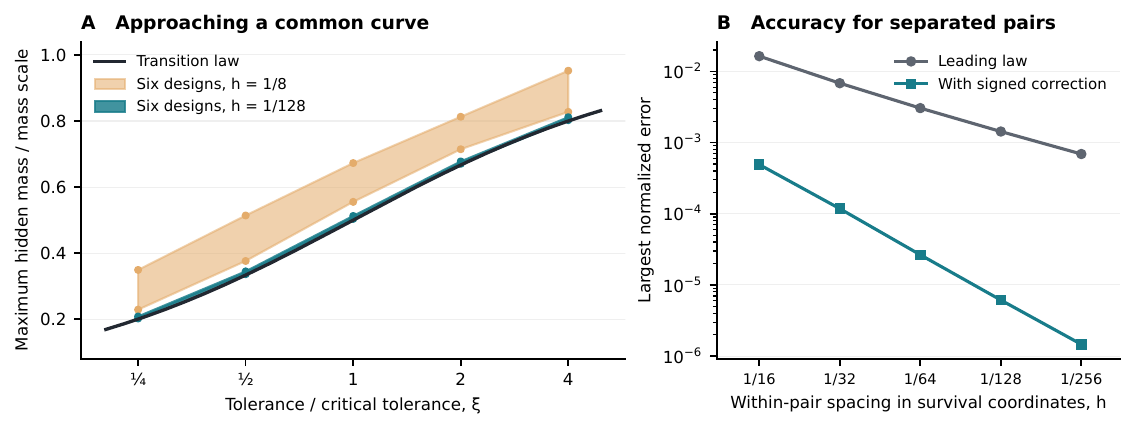}
\caption{How closely does the transition law predict a finite measurement? In A the line gives $\xi/(1+\xi)$. A vertical value of $0.5$ means that half the participating fraction $\omega$ can hide. The endpoints of each vertical range are the smallest and largest certified normalized masses among six designs at one tolerance. Shading connects the five tested tolerances. Reducing the physical-rate separation brings the points closer to the line. In B two pairs of survival coordinates remain separated as the spacing $h$ within each pair decreases. Each marker gives the largest absolute error in $A/\omega$ across four configurations and three tolerances. Including the signed term in Theorem~\ref{thm:multicenter} improves finite-parameter accuracy.}\label{fig:collapse8}
\end{figure}

Across the six designs, the largest absolute discrepancy in $A/\omega$ is below $0.181$, $0.0490$ and $0.0125$ at $h=1/8,1/32,1/128$, respectively. At $\xi=1$ the fluorescence-shaped design gives $0.5610$, $0.5166$ and $0.5042$, approaching the predicted half-mass level. These calculations scale down the observed \emph{relative} error profile. At $h=1/32$, for example, the required absolute half-width in the first channel is $4.82\times10^{-8}$. The measured high-count half-width is about $0.044$. If the measured half-widths are centered on a generated two-rate trace at their actual size, all nine fluorescence bands for equally weighted sites $X=0.84,0.87$ admit a single site at $X=0.855$. This remains true down to $11.243\%$ of the measured widths. The finite-band calculation detects saturation well outside the tested transition.

\paragraph{How much does the first correction buy?}
For two separated pairs we vary the distance between the pairs, the local weights, the channel precision profile and the tolerance. At within-pair spacing $1/256$ the largest normalized error in the leading law is $6.91\times10^{-4}$ to three significant digits. The largest error after the signed correction is below $1.47\times10^{-6}$ (Figure~\ref{fig:collapse8}B). All four-site maxima satisfy the finite-parameter vertex and positive-quadrature checks.

A second control centers synthetic echo readings on a known reference at $X=0.60$ and queries $X=0.69$. Its two half-widths, $0.01745$ and $0.01193$ to five decimals, come from the measured toluene echoes. A physical two-rate population and a quadratic upper certificate give the finite maximum in Table~\ref{tab:nmr-predictions}. The signed correction reduces the survival-coordinate error from $6.403$ to $0.252$ percentage points.

\begin{table}[htbp]\centering
\begin{tabular}{@{}lrr@{}}\toprule
Calculation & Signal at query (\%) & Absolute error (pp)\\\midrule
Certified finite maximum & 73.825 & 0.000\\
Physical-rate leading law & 82.375 & 8.550\\
Survival-coordinate leading law & 80.227 & 6.403\\
Survival law with signed correction & 74.077 & 0.252\\\bottomrule
\end{tabular}
\caption{A known reference at measured NMR precision. Errors are percentage points of total signal, calculated before rounding. The generated signal centers allow each forecast to be checked against the physical extremum.}\label{tab:nmr-predictions}
\end{table}

\paragraph{A finite-precision decision rule.}
Choose the acceptable error on the hidden mass before using a prediction $\omega p$. If a feasible population and a global certificate enclose the true maximum between $L$ and $U$, the error is at most $\max\{|L-\omega p|,|U-\omega p|\}$. Use the formula at the requested accuracy when this bound fits the tolerance. Otherwise report the certified interval or tighten it. For the NMR control, a tolerance of half a percentage point accepts the signed prediction. The proof of the asymptotic law supplies uniformity as the separation shrinks. This finite calculation supplies the numerical error at the chosen experiment.

\section{What the experiments resolve}\label{sec:illustrations}
The controlled calculations start from known reference rates. In a measured decay, the unknown quantity is the distribution itself. Two experiments illustrate the difference. A mixture of dyes has independent lifetime controls. Repeated NMR acquisitions measure instrument repeatability. A third example asks whether the positive-mixture model fits the observations at all.

\subsection{How precisely do the photons locate a known dye lifetime?}\label{sec:optical}
Kaye et al.\ measured a solution containing two dyes. Separate control measurements gave lifetimes of $0.48$ and $4.03$ ns \cite{kaye}. This gives the experiment a known long-lived component. Its share of the recorded light is unknown. How tightly can the mixed-solution trace locate that light?

The instrument recorded photon arrival times after excitation. Starting at $3.12710$ ns, beyond the last occupied bin of the measured pulse response, we counted photons in ten adjacent gates of width $\Delta=0.625419675$ ns. For an exponential with lifetime $\tau$, each gate retains the fraction $x=\exp(-\Delta/\tau)$ of the preceding gate's expected count. Thus, for a positive mixture, the ratio of gate $j$ to gate 0 is the mean of $x^j$ under the signal distribution in gate 0. A constant background has $x=1$. All percentages below refer to \emph{light in the first retained gate}. Relating light to dye concentration would require independent brightness and detection calibrations.

The archive supplies 300 random thinnings at each of two photon budgets from one measured acquisition. We use draw 1 at each budget for the comparison. For each of nine gate ratios we condition on the combined count in that gate and the reference gate. Hoeffding's inequality gives an interval for the conditional proportion, which we transform to a ratio. A union bound across both budgets gives at least $1-36e^{-7}>96.7\%$ simultaneous coverage for the designated draws. For finite-master subsampling the intervals target the master acquisition's empirical ratios. Population interpretation additionally requires a stationary arrival model and a representative instrument response.

The broad open interval $(1,8)$ ns contains the long control lifetime. How much of its signal can the readings establish on their own, and how much do the controls add? Table~\ref{tab:prior-optical} holds all nine error bands fixed within each column. With unrestricted lifetimes the higher-count designated histogram guarantees at least $54.49\%$ in the window. Across 300 thinnings at each budget, the medians of the certified lower bounds are $0\%$ and $42.12\%$, respectively.

\begin{table}[htbp]\centering
\begin{tabular}{@{}lrr@{}}\toprule
Allowed signal sources & 12,800 photons & 100,000 photons\\\midrule
Any positive lifetime mixture and background & $0\%$ & $54.49\%$\\
Two fixed control lifetimes and background & $63.8015\%$ & $85.4472\%$\\
Two fixed control lifetimes, zero background & $83.5132\%$ & $93.8602\%$\\\bottomrule
\end{tabular}
\caption{Minimum fraction of first-gate signal in $(1,8)$ ns, rounded to the precision shown. Each column uses draw 1. Reading across a row changes the photon budget. Reading down a column adds information about the allowed signal sources. With lifetimes fixed at 0.48 and 4.03 ns, only the 4.03 ns component contributes to this window.}\label{tab:prior-optical}
\end{table}

For the fixed-lifetime comparison, let $p$ be the long-lived fraction and $b$ the fraction of constant background. The predicted ratio in gate $j$ is $(1-p-b)e^{-j\Delta/0.48}+p e^{-j\Delta/4.03}+b$, with lifetimes in ns. Minimizing $p$ over the nine bands gives the middle row. Setting $b=0$ gives the last row and requires separate grounds for excluding background. These comparisons treat the control lifetimes as exact. As a sensitivity check, allowing each lifetime to vary independently by $\pm1\%$ gives minima of $63.47\%$ and $85.22\%$ with background, and $83.09\%$ and $93.28\%$ without it. The $1\%$ range is an illustrative perturbation. The reproduction materials give feasible weights and polynomial certificates for all four cases at both precisions.

\textbf{What happens if we narrow the window?} In the higher-count histogram, a population with $41\%$ of the signal at $x=0.730$ and $59\%$ at $x=0.902$ fits all nine error bands. Its lifetimes are $1.98728$ and $6.06375$ ns. Even the \emph{closed} window $[2,6]$ ns can therefore contain zero signal. At the other extreme, a single lifetime of $4.03$ ns also satisfies all nine bands at both budgets. Every open window containing the reference lifetime can hold $100\%$ of the signal. Both statements use the recorded response and the same error calibration.

To measure how closely the data locate the light, widen the open window $I_h=\{\tau:|\tau-4.03\,\mathrm{ns}|<h\}$. Intersect this window with positive finite lifetimes. Finite lifetime endpoints reverse order in the physical rate $k=1/\tau$; for example, $(1,8)$ ns corresponds to $(1/8,1)$ ns$^{-1}$. Finite windows exclude constant background. For each $h$ we bound the smallest compatible fraction from below by bounding the mass outside $I_h$ from above. An explicit feasible spectrum then shows how nearly this lower guarantee can be attained. At 100,000 photons the threshold half-width (the infimum of widths giving a $50\%$ guarantee) lies between $3.21875$ and $3.234375$ ns. At 12,800 photons it lies between $6.8984375$ and $6.921875$ ns. At the lower photon budget the half-width extends past zero lifetime, so the physical window is approximately $(0,10.94)$ ns. At the higher budget it is approximately $(0.80,7.26)$ ns. The lower-count bands even admit $26\%$ constant background, so every finite lifetime window can miss more than a quarter of the signal. This obstruction persists for arbitrarily long finite lifetimes.

\begin{figure}[htbp]\centering
\includegraphics[width=\linewidth]{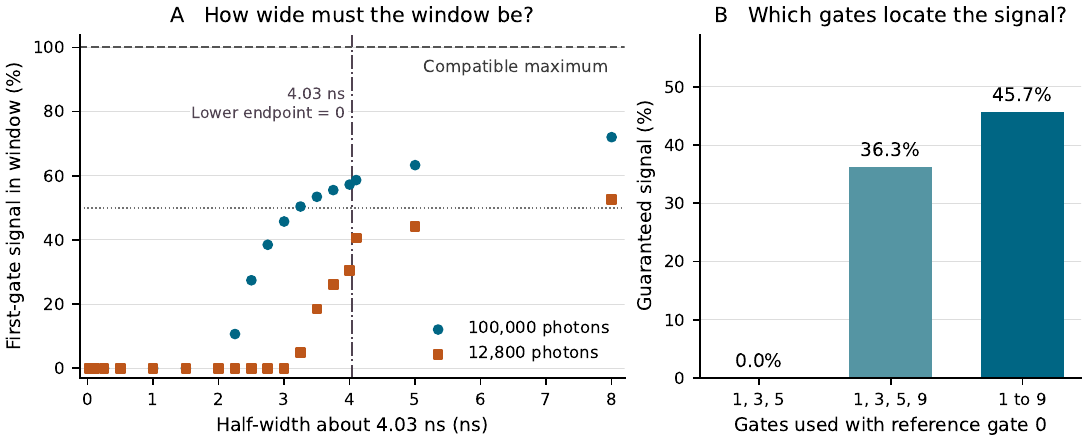}
\caption{How accurately do the measured photons locate the long-lifetime signal? A. Each symbol is a certified lower bound on first-gate light in an open window centered at $4.03$ ns, using all nine ratios in the designated histogram. The vertical line at $h=4.03$ ns marks where the lower endpoint reaches zero lifetime. The gap to an explicit feasible spectrum is below $0.002$ percentage points. The maximum is $100\%$ at both budgets for every positive half-width. B. For the window $(1.03,7.03)$ ns at 100,000 photons, the minimum compatible signal rises from zero to over $36.29\%$ on adding gate 9, and to over $45.74\%$ with all nine gates. Gate 9 spans $8.756$--$9.381$ ns after excitation. The comparison uses already recorded photons.}\label{fig:optical-window}
\end{figure}
This example separates three scientific questions. Detecting heterogeneity asks whether any single lifetime fits. Locating light asks how much signal every compatible spectrum places in a chosen range. Identifying a specific component requires enough information to assign a lifetime and amplitude to it. Here the calibrated counts locate light in a broad range, while independent dye controls supply the evidence for two components.

\begin{figure}[htbp]\centering
\includegraphics[width=.94\linewidth]{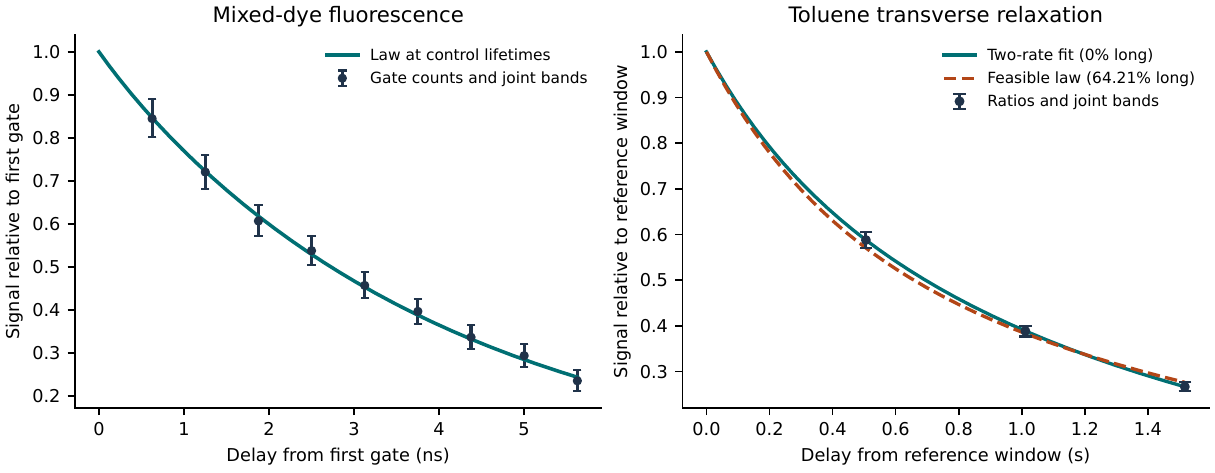}
\caption{Which decay traces can the instrument distinguish? Left, a mixture at the two independently measured lifetimes fits the higher-count fluorescence readings. Right, two populations fit the same three toluene ratios. One assigns zero signal to the slow region and the other assigns about 64.21\%. Symbols show measured ratios. Curves between observations are predictions.}\label{fig:measured}
\end{figure}

\subsection{What does one more NMR echo window reveal?}
For toluene, a third measured window reduces the largest admissible slow signal from $70.01\%$ to $64.21\%$. Here slow means $T_2\ge1.758$ s and the fraction refers to peak-voltage response in the reference window. Nonnegative least squares assigns $21.15\%$ of signal to this region, a Tikhonov fit assigns $21.28\%$, and a two-rate fit assigns zero. All three match the bands. The continuum certificate bounds every compatible positive spectrum (Figure~\ref{fig:measured}).

Huggins et al.\ recorded five technical repeats per hydrocarbon control \cite{huggins,arts}. Each reading averages 16 echoes. The reference starts at $0.05057$ s. The next three windows are translated by multiples of $\Delta=0.505689$ s, so $X=\exp(-\Delta/T_2)$ and the slow threshold is $X\ge3/4$. The measured ratios are $0.5880$, $0.3889$ and $0.2668$. Paired Student tests account for variation in the shared reference. A critical value of 6.5 with four degrees of freedom gives at least $95.37\%$ simultaneous coverage across four offsets in four controls under independent Gaussian repeat contrasts. The first two toluene bands already imply $\operatorname{Var}(X)\ge0.01034$. This establishes heterogeneity in the effective response. Identifying molecular environments requires independent physical evidence.

\emph{Receiver-floor sensitivity.} A window at $9\Delta$ averages about $6.15$ mV. Including its repeat band literally gives an upper bound below $9.94\%$. Allowing an additional $11.5$ mV of late-response uncertainty restores it to about $64.21\%$. A useful next control is a receiver-baseline measurement at the same gain and pulse settings. The third-window result uses earlier readings. The archive also tests all nine candidate lags, inter-channel bias and the other three controls.

\subsection{When do the observations call for a different model?}\label{sec:tablet}
Shah et al.\ measured cumulative release from twelve tablets at 30, 60, 90 and 180 minutes \cite{shah,dip2}. Keep each tablet's four readings together. If $u$ and $v$ are its remaining fractions at 60 and 90 minutes, every static positive mixture has $v\ge u^{3/2}\ge3u/4-1/16$. Each tablet violates the second inequality even after allowing half a percentage point of rounding at both times. The smallest surviving discrepancy is $0.275$ percentage points.

For population inference assume independent tablets and normal marginal errors. Four simultaneous Student intervals, including rounding in means and standard errors, have half-widths of $2.355$, $2.778$, $3.044$ and $2.586$ percentage points. The favorable band edges still violate $L_{30}^3\le U_{90}$. Restoring feasibility takes a common addition of $0.4167$ percentage points to all four half-widths. With boundedness and independence alone, a distribution-free 95\% calibration permits the full range $[0,1]$ for the two-tail fraction. The unique least-squares rate, $0.01573$ min$^{-1}$, has RMS error $2.41$ percentage points. The uncertainty and the observation model determine which interpretations the fit supports. Supplement, Section~\ref{app:tablet} documents the $3.97\%$ upper bound under an imposed 3.2-point tolerance. The supplement applies a related model check to a PET product yield \cite{pet}.

\section{From a fitted decay to a measurement decision}
In the dye mixture the readings require substantial emission within a broad lifetime band. Two finite lifetimes outside a narrower band reproduce the same observations. In toluene the first two windows establish heterogeneity in the effective decay and the third improves the upper bound on slow signal. In the tablet series each individual unit fails a necessary shape test. These outcomes call for different next steps: narrow a lifetime window, measure a receiver baseline, or revisit the kinetic model.

For an existing experiment, specify the signal normalization, a target with explicit endpoints, and a joint measurement-error band. Test feasibility, then compute examples and certificates for the requested fraction. For a proposed experiment with known reference components, $\omega$ gives the amount of signal at stake and $\eta_*$ gives a precision target. The dimensionless ratio $\eta/\eta_*$ predicts the fraction of $\omega$ that can hide as the rates approach one another. At finite spacing the benchmarks show how to check the prediction to a stated accuracy. Comparing future sampling times also requires an attainable precision and a cost for each time.

The transition laws cover a shrinking interior cluster on a fixed design and small coupling on a consecutive grid. A moving target window poses a further question. A single optimal distribution can put signal at several places in the window. The supplement gives an exact example and a computable condition under which a point certificate remains informative after widening the target. A sharp transition law for arbitrary moving windows on general grids remains open.

\section*{Computational methods and availability}
The reproduction materials contain the extraction procedures, interval conventions, reference inputs and certificate readers. Numerical optimization proposes distributions and coefficients. The readers check the original bounds and the dual inequalities across the continuous rate axis. The Lean~4 development gives machine-checked formulations of the three main results. Its theorem map specifies the definitions and assumptions at each entry point. The interpolation, compactness and quadrature arguments below give an independent mathematical derivation. Separate records document statistical calibration, algebraic examples and numerical replays.

\paragraph{Reuse.}
The author licenses the original article text and original figures under \href{https://creativecommons.org/licenses/by/4.0/}{CC BY 4.0} and the original code under the MIT License. Figure~\ref{fig:measured}, which includes adapted hydrocarbon echo data from Huggins et al.\ \cite{huggins,arts}, is licensed under \href{https://creativecommons.org/licenses/by-sa/4.0/}{CC BY-SA 4.0}. Source-specific terms and instructions for obtaining experimental records accompany the reproduction materials.

\section*{Acknowledgments}
\paragraph{Generative AI use.}
OpenAI Codex assisted with mathematical exploration, proof development, Lean formalization, numerical code, literature searches and writing. The original development recorded model identifiers \texttt{gpt-6.2-sol} and \texttt{gpt-6.2-astra}. The revision also used Codex to reorganize the proofs, check equation dependencies and compare the mathematical claims with earlier work. The software archive specifies the statements checked by the Lean kernel and the inputs used in numerical verification.

\paragraph{Data and code availability.}
The preprint is available as \href{https://arxiv.org/abs/2611.36950}{arXiv:2611.36950} and as Zenodo version 2.0, \href{https://doi.org/10.5281/zenodo.24757539}{10.5281/zenodo.24757539}. The \href{https://github.com/xarryzero/kinetic-rate-resolution}{code repository} provides \href{https://github.com/xarryzero/kinetic-rate-resolution/releases/tag/v1.1.0}{software release 1.1.0} with the formal proofs and reproduction materials. The release DOI is \href{https://doi.org/10.5281/zenodo.24757730}{10.5281/zenodo.24757730} and the all-versions software DOI is \href{https://doi.org/10.5281/zenodo.24567772}{10.5281/zenodo.24567772}.
\appendix
\section{Finite-error certificates and the physical endpoint}\label{app:finite}
\paragraph{Proof of Theorem~\ref{thm:face}.}
The affine hull of $M$ is $\mathbb R^n$ by the linear independence of the constant and the $n$ exponentials. A compactified optimizer exists by weak compactness and upper semicontinuity of the mass of a closed set.

Finite positive cubature gives at most $n+2$ atoms preserving the moments and mass of an optimizer \cite{didio}. Maximizing the reward again among its weights at fixed moments is a finite linear program with at most $n+1$ independent equalities. A basic optimum uses at most $n+1$ atoms. 

\paragraph{Interior contact and finite-rate repair.} The hypograph of the greatest reward at a fixed moment is closed convex by compactness and upper semicontinuity. When $B$ meets $\operatorname{int}M$, its relative interior does too. The finite-dimensional convex sum rule for the negative of this value and the box indicator gives an optimal affine majorant at an optimizing moment, with the appropriate signed band endpoints. Testing on point masses gives the dual in \eqref{eq:univ}. Interior duality and finite cubature are classical moment tools \cite{bp,didio}.

For an interior feasible moment $y_0$, a whole neighborhood belongs to the convex hull of a fixed finite collection of \emph{finite} rate vectors. Indeed the finite-rate convex hull is convex and dense in $M$. Choose an interior simplex around $y_0$ inside $M$, express its vertices by finite compactified mixtures and approximate any atom at infinity by a finite one. The perturbed simplex still contains a neighborhood. Now take a compactified optimizer at readings $y$. It can be made atomic while preserving all readings and target mass by finite positive cubature on $(1,F,{\bf1}_{T})$. Replace any mass at infinity by finite rates tending to infinity, choosing rates in $E$ if infinity is rewarded. This increases readings by a vector $e\to0$ and preserves or increases the target mass. For any $0<\theta<1$ sufficiently small errors permit a finite law with readings $y_0-(1-\theta)e/\theta$. Mix this repair law with weight $\theta$ and the replaced law with weight $1-\theta$. The resulting readings satisfy
\begin{equation}\label{eq:physical-repair}
 (1-\theta)(y+e)+\theta\{y_0-(1-\theta)e/\theta\}
 =(1-\theta)y+\theta y_0\in B .
\end{equation} A diagonal choice $\theta\downarrow0$ in \eqref{eq:physical-repair} proves equality of the suprema, even for zero-width coordinates. A physical optimizer, if present, can be made finite atomic by the same cubature principle.

\paragraph{Boundary contact and attainment of the dual.} Separate the open convex set $\operatorname{int}M$ from $B$ in the boundary case, obtaining $b\ne0$ with $c=\max_M b^\top y=\min_B b^\top y$. Each nonzero $b_i$ fixes the corresponding feasible reading to a particular endpoint of its interval. Write $h=c-b^\top F$, $Z=\{k:h(k)=0\}$, and let $m$ be the number of nonzero measured coefficients. Every feasible law is supported on $Z$. On finite real intervals a nonzero sum of $l$ exponentials has at most $l-1$ zeros counting multiplicity, by induction after dividing out one term and applying Rolle's theorem to the derivative. Interior zeros of a nonnegative contrast have even order. If infinity is a zero then $c=0$ and one fewer exponential remains. Thus $2|Z\cap(0,\infty)|+|Z\cap\{0,\infty\}|\le m$, including the possible zero at the finite endpoint. The evaluation matrix of the constant and $|Z|-1$ distinct active channels on $Z$ is nonsingular: a vanishing combination would violate the same zero count, dropping the constant if infinity is among the zeros. The active readings and total mass therefore fix every weight. Separation also certifies emptiness when the two compact convex sets are disjoint (clip bands to $[0,1]$ first).

At boundary contact use the exposing contrast $h\ge0$ just constructed, oriented so its maximal box expectation is zero. In $x=e^{-K}$ coordinates use only its active measured powers and the constant. At each positive zero of $h$, interpolate $V$ with value $\mathbf1_T$ and derivatives zero through one less than the order of the zero. If zero is a contact impose only its value there. These conditions have full row rank by the Chebyshev zero count: choose the constant and the requisite number of active powers; a zero at the origin removes the constant. Near a positive contact in $T$, $V-1=O(h)$. Near a contact outside $T$, closedness makes the indicator zero and $V=O(h)$. At an origin contact the lowest nonconstant power used by $V$ is at least the smallest active power of $h$, giving the same bound. Elsewhere away from contacts $h$ is bounded below. Thus $V+Ch$ majorizes the indicator for large $C$. Its measured coefficients then have the same signs as those of $h$. Their maximizing box endpoints are fixed by the unique feasible law, on whose support $V+Ch$ equals the reward. This proves dual attainment also at the boundary. For comparison, exact first and second moments of a point mass at an interior survival coordinate $a$ determine that law, yet the continuous reward $|x-a|$ cannot have a differentiable touching majorant at zero value. General continuous risks need not share indicator dual attainment on the boundary.

\section{Proofs of the transition laws}\label{app:cluster}
\noindent The first proof constructs a physical distribution that hides the predicted fraction at the query, then constructs an upper bound from the measured exponential channels. The second proof exploits consecutive observations to identify the exact extremizer and bound the signed correction. The quadratic estimate below is common to both arguments.
\subsection{A quadratic estimate shared by both proofs}\label{sec:energy}
Both results reduce the measured error to an extra direction in the same constrained quadratic problem. The normalized coordinates keep the leading energy finite when weights or spacings become very small. The following estimate gives the first correction without inverting a nearly singular raw moment matrix.
\begin{proposition}[Normalized trial energy]\label{lem:trial-energy}
Let $\|v\|=1$, $e=(v,1)$ and $G_0=\operatorname{diag}(I,\xi)$ with $\xi>0$.
On $e^\top u=1$, the minimum of $u^\top G_0u$ is $f(\xi)$ from \eqref{eq:master},
attained at $u_0=(\xi v,1)/(1+\xi)$. For a symmetric $G$, let $\Pi$ project
orthogonally onto $e^\perp$. Put $T=\Pi G|_{e^\perp}$ and $a=\Pi Gu_0$.
If $T\succeq cI$ for $c>0$, then
\begin{equation}\label{eq:trial-remainder}
 \min_{e^\top u=1}u^\top Gu=u_0^\top Gu_0-a^\top T^{-1}a,
 \qquad 0\le a^\top T^{-1}a\le\|a\|^2/c .
\end{equation}
\end{proposition}
\begin{proof}
The identity $G_0u_0=f(\xi)e$ gives stationarity for $G_0$. For $G$ write
$u=u_0+h$ with $h\perp e$. The excess over the trial value is
$2a^\top h+h^\top Th$. Completing the square gives the unique minimizing
displacement $-T^{-1}a$ and \eqref{eq:trial-remainder}.
\end{proof}
In particular, for $\xi$ in a compact subinterval of $(0,\infty)$, if
$G-G_0=\left(\begin{smallmatrix}B&k\\k^\top&0\end{smallmatrix}\right)$ with
$B=O(\epsilon^2)$ and $k=O(\epsilon)$, then $\|a\|=O(\epsilon)$ and $T$
stays uniformly positive. Evaluating the trial value in \eqref{eq:trial-remainder} gives
\begin{equation}\label{eq:trial-linear}
 \min_{e^\top u=1}u^\top Gu
 =f(\xi)+\frac{2\xi}{(1+\xi)^2}v^\top k+O(\epsilon^2).
\end{equation}
The unsigned bound in \eqref{eq:trial-remainder} is relative to the full trial energy,
which includes its upper-left-block contribution.

\begin{lemma}[Quadrature with a prescribed node]\label{lem:anchored}
Let $L$ be a real linear functional on real polynomials of degree at most $2r$ with
$L(1)=1$ and positive definite Gram form on degree at most $r$.
Let $q$ minimize $L(q^2)$ subject to $q(z)=1$, and suppose $\deg q=r$.
Write $Q$ for its monic normalization. Then $z$ and the roots of $Q$ are
$r+1$ distinct real nodes of a positive quadrature exact through degree $2r$.
Its weight at $z$ is $\rho=L(q^2)$. If $P=\prod_i(X-x_i)$ has distinct real
zeros and $Q=P+\sum_i a_i P/(X-x_i)$, then each zero $y$ of $Q$ satisfies
\begin{equation}\label{eq:root-localization}
 \min_i|y-x_i|\le\sum_i|a_i|.
\end{equation}
\end{lemma}
\begin{proof}
Stationarity gives $L(qp)=\rho p(z)$ for $\deg p\le r$ by subtracting $p(z)q$.
Thus $F=(X-z)Q$ satisfies $L(Fs)=0$ for $\deg s\le r-1$.
A repeated real root or a conjugate pair would give $F=g s$ for a monic quadratic
$g=(X-b)^2+c^2$, $c\ge0$. Then $L(Fs)=L(((X-b)s)^2)+c^2 L(s^2)>0$,
a contradiction. Divide a polynomial of degree at most $2r$ by $F$ and
interpolate the remainder at its simple real roots. The weights are
$L(l_y)=L(l_y^2)>0$ since $l_y^2-l_y$ is divisible by $F$ for each cardinal
polynomial $l_y$. At $z$, $l_z=q$. Finally $1\le\sum_i|a_i|/|y-x_i|$
at a nonreference root by division by $P(y)$, proving \eqref{eq:root-localization}.
\end{proof}
The full-degree hypothesis is essential. For equal weights at $1/4$ and $3/4$ and anchor $z=1/2$, the Gram form for $r=1$ is positive definite but the minimizer is the constant one. Any positive two-site law containing its own mean has zero variance, so it cannot realize these moments.
For fixed moments the prescribed-node construction and reciprocal Christoffel
weight are classical \cite[Theorem 3.1]{hurlimann}; see also \cite[Theorem 10.21]{schmudgen} and \cite{simon}. We use
\eqref{eq:root-localization} to keep the nodes physical while the moment
functional itself varies with the error box.

\subsection{A cluster of physical rates}
Use the notation of Section~\ref{cl:section} and abbreviate $M=2r$, $\Gamma=\eta/d^2$, so $\xi=\Gamma C$.
\begin{proof}[Proof of Theorem~\ref{cl:nested}]
Use centered physical displacement $x$, $\phi_j(x)=g_j e^{-t_jx}$, and $L_0 f=\sum_i w_i f(x_i)$. The linear program
\begin{equation}\label{eq:design-lp}
 \max u_M:\quad \left|g_j\sum_{m=1}^{M}(-t_j)^m u_m/m!\right|\le\sigma_j,\quad 1\le j\le n
\end{equation}
has value $C$ in \eqref{eq:design}. The dual minimizes $\sum_j\sigma_j|c_j|$ subject to the corresponding moment equalities for coefficients of degrees 1 to $M$. A basic solution uses at most $M$ nonzero coefficients and fewer cannot annihilate degrees $1,\ldots,M-1$ by Vandermonde invertibility. For each $M$-subset the coefficients are the top Lagrange divided-difference row, including time zero for normalization. This yields \eqref{eq:design} from \eqref{eq:design-lp}. Fix one maximizing vector $u^*$.

For $\theta\in[1/2,1]$, put $u=\theta u^*$, and define a perturbed functional $L$ on polynomials through degree $M$ by $L(1)=1$, $L(x^m)=L_0(x^m)+\eta u_m$ for $m>0$. In the degree-$r$ basis
\begin{equation}\label{eq:weighted-basis}
 \ell_1/\sqrt{w_1},\ldots,\ell_r/\sqrt{w_r},\ P/d ,
\end{equation}
the evaluation vector at $z$ is $\sqrt S\,e$, where $e=(\ell_1(z)/\sqrt{w_1S},\ldots,\ell_r(z)/\sqrt{w_r S},1)$ has norm $\sqrt2$. The Gram matrix of $L$ is
\begin{equation}\label{eq:cluster-gram}
 G=\operatorname{diag}\big(1,\ldots,1,\eta u_M/d^2\big)+O(H)
\end{equation}
in operator norm. Indeed, $|d/(\sqrt{w_i}P'(x_i))|=|z-x_i|\,|\ell_i(z)|/\sqrt{w_i S}\le 2H$. After multiplying each of the first $r$ basis elements of \eqref{eq:weighted-basis} by $d$, its coefficients are therefore $O(H)$; the coefficients of the monic polynomial $P$ are bounded, and all but its top coefficient are $O(H)$. In the last diagonal Gram entry the perturbation of $P^2$ hence gives $\eta(u_M+O(H))$. The unperturbed Gram form has diagonal $(1,\ldots,1,0)$. All estimates stay uniform if a reference weight tends to zero at an arbitrary speed.

Since $\eta u_M/d^2=\theta\xi$, \eqref{eq:cluster-gram} is uniformly positive definite. Proposition~\ref{lem:trial-energy} identifies the leading constrained energy in this normalized basis. The least $L(q^2)$ subject to $q(z)=\sqrt S$ is $\rho_0=(e^\top G^{-1}e)^{-1}=(1+d^2/(\eta u_M))^{-1}+O(H)$. Thus the least square cost at value one is $\rho=\rho_0/S$. The coefficients of the minimizer $q$ normalized at $\sqrt S$ are $\rho_0 G^{-1}e$. They are bounded and their last entry $b$ is bounded away from zero, since the last entry of $G^{-1}e$ is $d^2/(\eta u_M)+O(H)$ and $\rho_0$ is bounded below. The monic normalization $Q=(d/b)q$ satisfies
\begin{equation}\label{eq:kernel-monic}
 Q(x)=P(x)+\sum_i O(H)\prod_{j\ne i}(x-x_j),\qquad Q(z)=D/b .
\end{equation}
By \eqref{eq:kernel-monic} and Lemma~\ref{lem:anchored}, $L$ has a positive
quadrature $\mu$ with weight $\rho$ at $z$. Each of its other nodes is within
$O(H)$ of a reference displacement by \eqref{eq:root-localization}. Hence
all rates are finite and positive for small $H$, with no lower bound on a
relative gap or a reference weight.

To check the actual readings, interpolate $\phi_j$ by a polynomial $p_j$ of degree $M-1$ at each $x_i$ doubly. On the supports of both laws the Hermite remainder is
\begin{equation}\label{eq:hermite-transfer}
 \phi_j(x)-p_j(x)=P(x)^2\left(\phi_j^{(M)}(0)/M!+O(H)\right).
\end{equation}
Uniformity, including arbitrarily close nodes, follows from the divided-difference mean-value formula \cite{deboor}. The coefficients of $p_j$ are $\phi_j^{(m)}(0)/m!+O(H)$, $0\le m<M$, by its Newton divided-difference representation. The remainder is zero on the reference support and $L(P^2)=\eta(u_M+O(H))$. Integrating \eqref{eq:hermite-transfer} therefore gives
\begin{equation}\label{eq:reading-transfer}
 \int\phi_j\,d\mu-L_0\phi_j=
 \eta\sum_{m=1}^{M}\frac{\phi_j^{(m)}(0)}{m!}u_m+O(\eta H).
\end{equation}
The remainder in \eqref{eq:reading-transfer} is proportional to the allowed error even when $D^2$ is much smaller than $H^{2r}$. Choosing $\theta=1-KH$ with $K$ larger than the uniform remainder constant divided by $\min_j\sigma_j$ enforces all bands. This gives the lower bound.

For the upper bound use the uncontracted vector and the resulting $Q$. Fix a minimizing subset $I$ and write the roots of $Q$ as $a_1,\ldots,a_r$. In the span of $1,\{\phi_j:j\in I\}$ choose $\Psi$ to vanish doubly at the roots and to have $\Psi^{(M)}(0)/M!=1$. This system has a unique bounded solution for small $H$. To verify uniformity, impose the first $M$ homogeneous equations as consecutive Newton divided differences of orders $0,\ldots,M-1$ on the list $a_1,a_1,\ldots,a_r,a_r$. Each row differs by $O(H)$ from the derivative row at zero divided by its factorial. The limiting augmented matrix is an invertible exponential Vandermonde matrix. Its last inverse column gives the coefficients with weighted absolute sum $C$ over measured rows, so the finite sum is $C+O(H)$. The argument of \eqref{eq:hermite-transfer} now gives $\Psi(x)=Q(x)^2(1+O(H))$ uniformly at the reference points and anchor. In particular $\Psi(z)>0$. A nonzero linear combination of $M+1$ exponentials, including the constant, has at most $M$ real zeros counted with multiplicity, by repeated Rolle differentiation after division by a nonvanishing exponential. Its $r$ double zeros exhaust the count. Thus $\Psi$ is globally nonnegative, also by continuity at an infinite rate. In particular the quotient is a certificate on the entire compactified axis:
\begin{equation}\label{eq:global-certificate}
 \frac{\Psi(x)}{\Psi(z)}\ge {\bf1}_{\{z\}}(x),
 \qquad -k_*\le x\le\infty .
\end{equation}

Integrating \eqref{eq:global-certificate} as in \eqref{eq:univ} bounds $A$ by $\{(1+O(H))L_0 Q^2+\eta(C+O(H))\}/\{Q(z)^2(1+O(H))\}$. Since each root of $Q$ is $O(H)$, $L(Q^2)=L_0 Q^2+\eta(C+O(H))$. Moreover $L(Q^2)/Q(z)^2=\rho$, while $\eta/Q(z)^2=\Gamma b^2/S=O(1/S)$ and $\rho=O(1/S)$. This gives $A\le\rho+O(H/S)$ and proves \eqref{cl:formula}.
\end{proof}

\begin{corollary}[Exact finite vertex on a controlling subgrid]\label{cl:corner}
Under Theorem~\ref{cl:nested}, put $M=2r$. If $n=M$ then for all sufficiently small $H$ the exact optimizer on the half-line is unique, has $r+1$ atoms, and reaches the lower sides at odd-numbered times and the upper sides at even-numbered times in increasing order. With more rows, the same conclusion on a minimizing $M$-subset applies if the unique moment perturbation saturating its alternating vertex lies strictly inside every unused linearized band.
\end{corollary}
\begin{proof}
Let $B_{jm}=\phi_j^{(m)}(0)/m!$ on the selected square subset and $v_j=\sigma_j\operatorname{sign}(c_j)$ for its nonzero top divided-difference coefficients $c_j$. Then $u^*=B^{-1}v$ has $u^*_M=C$. For $u$ in a fixed small ball around $u^*$, the Gram estimate \eqref{eq:cluster-gram}, Lemma~\ref{lem:anchored} and \eqref{eq:reading-transfer} hold uniformly, since $u_M$ stays positive. For each fixed configuration the quadrature depends continuously on $u$: its distinct simple roots vary continuously, and its weights solve an invertible Vandermonde system. Thus its vector of raw deviations divided by $\eta$, denoted $T(u)$, is continuous and satisfies $T(u)=Bu+O(H)$ uniformly. The map $u\mapsto u+B^{-1}(v-T(u))$ takes the closed ball of radius $K H$ around $u^*$ into itself for fixed sufficiently large $K$. Brouwer's theorem gives an exact solution $T(u)=v$. The other bands, when present, follow from their strict slack.

Use the exponential certificate $\Psi$ interpolated at the $r$ roots for this very quadrature. Its coefficients on the selected readings have the limiting signs of $c_j$, all nonzero. Hence its integral is largest at the alternating vertex. The constructed quadrature attains this bound. Equality forces this same vertex and support only at the anchor and the $r$ roots, since $\Psi$ is strictly positive off its specified zeros. Their weights are unique by the Chebyshev property. The sign pattern follows from $M$ even and the order of the selected times.
\end{proof}

\subsection{Consecutive readings and separated centers}
In the following proof $\Gamma=\eta/d^2$, $\alpha=\Gamma C=\xi$, $M=M_1$ bounds the widths, and $|d|/(\sqrt{w_i}|P'(x_i)|)=|z-x_i|\sqrt{\pi_i}$.
\begin{proof}[Proof of Theorem~\ref{thm:multicenter}]
Use the vertex \eqref{eq:alternating-vertex}. Let $L_*(X^j)=m_j^*$ and let $G$ be the Gram matrix of $L_*(fg)$ in the basis
$b=(\ell_i/\sqrt{w_i},P/d)$. On the reference law the matrix is
$\operatorname{diag}(I_r,0)$. The perturbation acts on $f$ as
$\eta\sum_{j\ge1}(-1)^j\sigma_j[X^j]f$. The coefficients of $R_i$ and $P$
are bounded in terms of $r$, and all coefficients of $P^2$ alternate in sign.
The perturbation of the upper-left block is therefore $O(\epsilon^2)$,
its off-diagonal column is $O(\epsilon)$ and its bottom-right entry is exactly $\Gamma C$.
Here $c\le C\le M4^r$, which gives $G\succ0$ and
$G^{-1}=\operatorname{diag}(I_r,1/(\Gamma C))+O(\epsilon)$.
The first $r$ entries of $u=b(z)/\sqrt S$ have norm one, and its last entry is one.
Thus $u^\top G^{-1}u=1+1/(\Gamma C)+O(\epsilon)$.
For the correction write $u=(v,1)$, $G=\left(\begin{smallmatrix}I+B&k\\k^\top&\alpha\end{smallmatrix}\right)$.
Equation~\eqref{eq:trial-linear} gives the normalized constrained minimum $f(\alpha)+2\alpha v^\top k/(1+\alpha)^2+O(\epsilon^2)$. The degree-$2r-1$ polynomial $R_iP$
has coefficient signs $(-1)^{j+1}$. Thus
$k_i=-\Gamma d C_i/(\sqrt{w_i}P'(x_i))$, $v_i=d/((z-x_i)\sqrt{w_i}P'(x_i))$,
so $\pi_i=v_i^2$ and
\begin{equation}\label{eq:signed-coupling}
 v^\top k=-\Gamma\sum_i\pi_i(z-x_i)C_i .
\end{equation}
Substituting \eqref{eq:signed-coupling} into \eqref{eq:trial-linear} with $\Gamma=\alpha/C$ gives \eqref{eq:correction} once attainment is established.

It remains to justify that the vertex is feasible and optimal on the physical interval.
The last entry of $G^{-1}u$ is bounded away from zero. The polynomial
$b(X)^\top G^{-1}u$ can therefore be rescaled to a monic polynomial
$Q=P+\sum_i O(\epsilon)R_i$. Lemma~\ref{lem:anchored} now gives positive quadrature at $z$ and the $r$
roots of $Q$. By \eqref{eq:root-localization} each residual node is within
$O(\epsilon)$ of a reference site in $[a,b]\subset(0,1)$, so it belongs
to the physical interval. The weight at $z$ is the minimum square cost,
$1/(S u^\top G^{-1}u)$, which becomes \eqref{eq:inverse-gram} in the
monomial basis.

Since all roots of $Q$ are positive, the coefficients of $Q^2$ strictly alternate.
This sign test concerns this particular square: even for $r=1$ the difference
between the corner and reference Gram matrices has determinant
$-(\eta\sigma_1)^2$, so a Loewner comparison of the two matrices fails
when $\eta\sigma_1>0$. The coefficient signs imply that the vertex maximizes the expectation of $(Q/Q(z))^2$ over the entire box.
This square majorizes the point indicator, and the quadrature attains its bound.
Equality forces support on the quadrature nodes and forces every nondegenerate
band endpoint; interpolation fixes all the weights.
\end{proof}

\paragraph{Exact remainder for two readings.}
For $r=1$ write $x$ for the reference survival, $s=z-x\ne0$,
$e_j=\eta\sigma_j$, $\alpha=(e_2+2xe_1)/s^2=\xi$ and $\beta=e_1/s$.
When $|\beta|\le1/4$, $\beta^2<\alpha$, and $z$ and
$y=x-s(\alpha+\beta)/(1+\beta)$ lie in $(0,1]$, the optimal mass is
$p=(\alpha-\beta^2)/(1+\alpha+2\beta)$. Indeed,
$p\delta_z+(1-p)\delta_y$ has moments $x-e_1,x^2+e_2$ and meets the
upper certificate $(X-y)^2/(z-y)^2$. Substitution gives the exact remainder
\begin{equation}\label{eq:two-reading-remainder}
 p-\left[f(\alpha)-\frac{2\alpha\beta}{(1+\alpha)^2}\right]
 =-\frac{\beta^2(\alpha-1)^2}{(1+\alpha)^2(1+\alpha+2\beta)} .
\end{equation}
By \eqref{eq:two-reading-remainder}, the linear approximation bounds the exact mass from above, with
error at most $2\beta^2$. At $\alpha=1$ it is exact throughout this regime.
Under the compact design bounds of Theorem~\ref{thm:multicenter},
$|\beta|\le |s|\,(\sup\xi)M_1/c$ and both physical nodes lie in $(0,1)$
for sufficiently small $|s|$, uniformly in the admissible widths.

\begin{proposition}[Agreement under the survival-coordinate transformation]
\label{lem:coordinates}
Take $2r$ consecutive times $j\tau$ in the interior single-cluster regime, with fixed positive
widths $\sigma_j$. Let $a=e^{-\tau k_*}$ and $g=\tau a$. Subscripts $K$ and $X$ distinguish
quantities in rate and survival coordinates. Then
\begin{equation}\label{eq:coordinate-transfer}
 S_X=S_K(1+O(H)),\qquad
 d_X^2=g^{2r}d_K^2(1+O(H)),\qquad
 \frac{\eta C_X}{d_X^2}=\frac{\eta C_K}{d_K^2}(1+O(H)).
\end{equation}
The constants are uniform over all relative gap sizes and positive weights.\end{proposition}
\begin{proof}
Each nonzero difference in $X=e^{-\tau K}$ is $-g$ times the corresponding difference
in $K$, up to relative error $O(H)$ by the mean-value theorem.
Substituting into each Lagrange ratio and each positive term of the harmonic formula
for $d^{-2}$ proves the first two identities. The consecutive divided-difference row gives
$C_K=\tau^{-2r}\sum_{j=1}^{2r}\sigma_j\binom{2r}{j}a^{-j}$.
In survival coordinates the monic annihilating polynomial has roots $a+O(H)$,
so $C_X=\sum_{j=1}^{2r}\sigma_j\binom{2r}{j}a^{2r-j}+O(H)=g^{2r}C_K+O(H)$.
Positivity of $C_K$ gives the last identity in \eqref{eq:coordinate-transfer}.
Finally, if $z_X$ and $x_{i,X}$ denote the survival anchor and reference nodes, then
$\epsilon=\max_i |z_X-x_{i,X}|\sqrt{\pi_i}\le\max_i|z_X-x_{i,X}|=O(H)$.
Thus both normalized transition predictions agree to the accuracy of the physical-rate law.\end{proof}

\clearpage
\begingroup
\setcounter{section}{0}
\setcounter{equation}{0}
\setcounter{figure}{0}
\setcounter{table}{0}
\renewcommand*{\thesection}{S\arabic{section}}
\renewcommand*{\theequation}{S\arabic{equation}}
\renewcommand*{\thefigure}{S\arabic{figure}}
\renewcommand*{\thetable}{S\arabic{table}}
\renewcommand*{\theHsection}{supplement.\arabic{section}}
\renewcommand*{\theHequation}{supplement.\arabic{equation}}
\renewcommand*{\theHfigure}{supplement.\arabic{figure}}
\renewcommand*{\theHtable}{supplement.\arabic{table}}
\renewcommand*{\thesupplementtheorem}{S\arabic{supplementtheorem}}
\let\theorem\supplementtheorem
\let\endtheorem\endsupplementtheorem
\section*{Mathematical supplement}

\noindent The support-count proof and computational certificates below accompany
Theorem~\ref{thm:face}. Subsequent sections treat window and prescribed-point
extrema, design comparisons and additional scale regimes.
\section{A sharp support count for unions of intervals}\label{irregular:section}
Write $x=e^{-K}$ in fixed time units and $I=[l,1]$, $0\le l<1$.
Positive $l$ imposes a physical upper rate, while $l=0$ includes the compactifying point.
The $n\ge1$ distinct positive reading times $t_i$ give $\phi(x)=(x^{t_1},\ldots,x^{t_n})$.
For a union $A\subset I$ of $b\ge1$ disjoint nonempty closed intervals (singletons allowed),
maximize $\mu(A)$ over probability laws with $\int\phi\,d\mu$ in a feasible closed box $B$.
Upper semicontinuity and compactness give attainment. A certificate $Q=c_0+\sum c_i x^{t_i}$
majorizes $\mathbf1_A$. Its upper bound equals $c_0+\sup_{y\in B}\sum c_i y_i$.
Equality forces support on the contact set and fixes the observed moment at the appropriate
band endpoint for every nonzero $c_i$.

\begin{theorem}[Sharp number of rates for interval unions]\label{irregular:support}
For every specified grid of $n$ times and each $I$ above, the worst-case minimum number
of atoms in an optimizing law, over feasible boxes and targets with $b$ components, is
\[
 S(n,b)=\min\{n+1,\ \lfloor n/2\rfloor+2b+1\}.
\]
The worst case can be realized with strictly positive widths on all bands and with all $b$
intervals nondegenerate. For a nonconstant optimal certificate with $m$ nonconstant terms,
its $k$ contact points satisfy $2k\le m+e+T$. Here $e$ counts contacts at a domain endpoint
or target boundary and $T$ counts changes of the indicator among the ordered contacts.
\end{theorem}
\begin{proof}
A nonzero combination of $s$ distinct powers has at most $s-1$ positive zeros counted with
multiplicity: divide by the smallest power, differentiate and apply Rolle's theorem
inductively. For a system including the constant, a zero at $0$ may be counted once by
deleting the vanishing constant coefficient before applying the positive zero bound.
This also implies nonsingularity of Hermite interpolation with interior derivative
conditions. The derivative of $Q$ has at most $m-1$ positive zeros. All contacts are among
the zeros of $Q$ or $Q-1$, hence finite. Each of the $k-e$ ordinary interior contacts is a
stationary point. An additional derivative zero lies strictly between consecutive contacts
of the same height. Thus $(k-e)+(k-1-T)\le m-1$. Since $e\le2b+2$, $T\le2b$, this
proves the bound whenever a nonconstant certificate is available. Carath\'eodory reduces
the support to $n+1$ while preserving every moment within the contact set.

Write $M=\operatorname{conv}\phi(I)$. It has nonempty interior by independence.
If $B$ meets $\operatorname{int}M$ and the maximum is below one, a nonconstant certificate
exists. Indeed the upper hypograph of attainable moment--mass pairs is closed convex. Its
upper boundary as a function of $y\in M$ is concave and upper semicontinuous. The convex
sum rule applied to its negative and the box indicator gives a supporting affine majorant at
the optimizing $y$ with the required endpoint signs: $\operatorname{ri} B$ intersects
$\operatorname{int}M$ by taking a short segment from $B\cap\operatorname{int}M$ toward
$\operatorname{ri}B$. Testing this majorant on point masses gives the certificate.
A constant majorant for a nonempty $A$ cannot attain a value below one.
If $B$ only touches $\partial M$, separation gives a nonconstant nonnegative generalized
polynomial vanishing on the support of every feasible law: subtract the separating
functional's support value and choose its sign. Interior zeros count at least twice,
endpoints at least once. There are at most $\lfloor(n+2)/2\rfloor$ distinct zeros.

If the maximum is one, fix a feasible exact moment $y$ generated entirely on $A$.
If $A$ consists of isolated points, ordinary Carath\'eodory gives at most $\min(n+1,b)$
atoms. Otherwise $M_A=\operatorname{conv}\phi(A)$ has nonempty interior. At its boundary,
a nonzero supporting polynomial is nonnegative on $A$ and zero at all representing nodes.
For $y$ in its interior, maximize a new channel $x^\tau$ at that fixed $y$ with $\tau>0$
distinct from all the times. Its finite concave value has a supporting affine majorant at
the interior point. Subtract $x^\tau$ from that majorant on $\phi(A)$ to obtain a nonzero
generalized polynomial nonnegative on $A$, zero on an optimizing representation.
There are at most $n+1$ zeros counted with multiplicity. All except at most $2b$ endpoints
count at least twice. Thus some optimizer uses at most $\lfloor(n+1+2b)/2\rfloor$ atoms,
also subject to the $n+1$ bound.

For sharpness first use $2d$ times, $d\ge1$, and put $b'=\min(b,d)$.
Choose $d+1$ ordered valley nodes, including both endpoints of $I$. Prescribe value zero
at both endpoints and at $b'-1$ interior valleys, and value one at the other $d-b'$ valleys.
Prescribe zero first derivative at every interior valley. These $2d$ Hermite conditions on
the $2d+1$ dimensional space including the constant are independent: the first $2d$ basis
functions have an invertible interpolation matrix by the zero count. The homogeneous solution
$h$ is unique up to scaling. Its value vanishes at the endpoints, its interior valley zeros
are exactly double, and it has no further zeros. Choose its sign positive between valleys.
To any interpolant add a sufficiently large positive multiple of $h$ that the value at one
selected point of each inter-valley gap exceeds one. Each gap now has an interior
maximum. These $d$ maxima and $d-1$ stationary valleys use up the entire derivative zero
count. Each intervening arc is strictly monotone and $Q$ is nonnegative.
The superlevel set $A_0=\{Q\ge1\}$ consists of $b'$ nondegenerate intervals. There are
$b'+1$ zero valleys, $d-b'$ unit valleys and $2b'$ level crossings. These give the full
contact list of $k=d+2b'+1$ points. The derivative count forces all $2d$ nonconstant
coefficients to be nonzero. If $b>b'$, delete $b-b'$ small, disjoint open subintervals
strictly within a region where $Q>1$. The resulting target has $b$ intervals and keeps
the same contacts.

For $n=2d$ or $2d+1$, apply the construction with any $2d$ of the specified times.
On the contact list take all $n$ moments, including the remaining time in the odd case.
Any $\min(k,n+1)$ distinct images are affinely independent by the generalized Vandermonde
property. Their full convex hull has dimension $s-1$ for $s=\min(k,n+1)$. Choose in its
relative interior a vector $y$ outside the finite union of hulls generated by fewer than
$s$ contacts. It has positive distance from that compact union and admits a representation
positive on every contact (subtract a small multiple of the full barycenter and renormalize).
Set each used row's upper endpoint to $y_i$ if $c_i>0$, and its lower endpoint to $y_i$
if $c_i<0$, with positive width on the other side. For the unused time, if any, give it a
positive-width band centered at $y_i$ with half-width less than the preceding distance.
Equality with the integrated certificate holds for the chosen law. It forces every optimum
to use only contacts and to match $y_i$ on all $2d$ used rows. Every permitted optimal
moment vector is thus too close to $y$ to be representable by fewer than $s$ contacts.
Since $s=S(n,b)$, this proves sharpness. For even $n$ with $k\le n+1$, all the weights
are unique by affine independence.
For $n=1$, take a nondegenerate closed target interval below the upper domain endpoint in
$x$. Set the reading's lower bound strictly between the channel values at the interval's
right end and the upper domain endpoint. The affine function in $x^{t_1}$ equal to one
at the former and zero at the latter is a majorant. A law giving positive mass to both
contacts to match the lower band endpoint reaches its bound. Additional intervals strictly
to the left, disjoint and avoiding contacts, preserve the example. This gives $S(1,b)=2$.
\end{proof}
This is a worst-case minimum over optimizing laws. Each finite physical cap has the same
sharp result. In the compactification an optimizer putting mass at $x=0$ requires a
separate attainment analysis in the space of finite rates.

In particular this gives an exact finite-dimensional primal formula without first
finding the supporting hyperplane. Write $s=S(n,b)$. For each integer $j=0,\ldots,s$,
maximize $\sum_{i=1}^j w_i$ over the simplex of $s$ weights with
$x_1,\ldots,x_j\in A$, $x_{j+1},\ldots,x_s\in I$, and
$\sum_i w_i\phi(x_i)\in B$. Take the largest value among the feasible programs.
Every objective so obtained is bounded above by the target mass of its own law.
Conversely one can order a sharp atomic optimizer with all target nodes first
and pad with zero weights. Each problem is compact, has at most $2s-1$ real
parameters after normalization, and can be nonconvex. The sharpness concerns the
atom count. It makes no claim about optimal arithmetic running time.

\paragraph{Checking a certificate for an arbitrary compact target.}
Nonnegativity on $I$ requires testing both endpoints and the interior critical points. On
each closed monotonicity interval of $Q$ test $Q\ge1$ at the downhill-most point of $A$
in that interval, if the intersection is nonempty. These at most $m$ target tests suffice
for any compact $A$. Finite interval unions allow direct intersection queries; an abstract
closed set requires an extremal-point oracle. Critical points can be isolated in exact real
arithmetic by repeatedly dividing by the lowest power, differentiating, and solving at most
one strictly monotone equation per resulting interval, checking stationary endpoints for
multiple roots. Exact sign or equality information is needed for coincidences. The
sign-variation version of the zero bound can further reduce the count
\cite{jameson}.

\section{Computational certificates}
\subsection{The nested example at finite scale}\label{app:nested-examples}
\paragraph{A finite remainder certified on a full interval of scales.}
For the example at $\eta=h^{10}/728$ the exact finite answer satisfies
\begin{equation}\label{cl:intervalrate}
 2/5<A_h<2/5+h\qquad (0<h\le1/32).
\end{equation}
Here is a reproducible integer-polynomial certificate. Write $u=1/2$, $b=u/(1+h^2)$, $X=e^{-K}$, and let $m_0=1$ and
\[
 m_j=\frac{u^j}{3}\{1+(1+h^2)^{-2j}+(1+h)^{-j}\}+(-1)^j h^{10}/728,\quad 1\le j\le6.
\]
These are the moments at the alternating raw vertex. Set $G=(m_{i+j})_{0\le i,j\le3}$ and $v=(1,b,b^2,b^3)^\top$. The exact value is $\rho=1/(v^\top G^{-1}v)$. Indeed all leading minors of $G$ are positive. If $q$ is the polynomial with coefficients $\operatorname{adj}(G)v$, every coefficient of $q(-X)$ is positive and every coefficient of $q(1+X)$ is negative. Consequently $q$ has degree three, and its three roots, real and simple by the preceding rank-one quadrature construction, all lie in $(0,1)$. The square of $q/q(b)$ bounds the anchor indicator, its nonconstant coefficients have alternating signs, and its integral over the whole raw box is largest at this vertex. The residual quadrature gives equality.

For mechanical verification, clear the common moment denominator $L=3\cdot64\cdot728(1+h^2)^{12}(1+h)^6$, put $W=8(1+h^2)^3$, and form
\[
 N=W^2\det(LG),\qquad E=L(Wv)^\top\operatorname{adj}(LG)(Wv).
\]
Then $\rho=N/E$. The code verifies all minor and coefficient signs, and positivity of $E$, $5N-2E$, and $(2+5h)E-5N$. For each integer polynomial it removes its initial power of $h$. If the remaining degree is $n$, the substitution $h=t/[32(1+t)]$, followed by multiplication by $(32(1+t))^n$, yields coefficients all of the required strict sign; the leading coefficient covers the right endpoint. Only integer sums and products are used. An independent nine-point check builds the alternating moment form as the reference measure plus $\eta(\delta_{-1}-\delta_0)$, evaluates the Lagrange kernel after adding the node $-1$, then applies a rank-one downdate at zero. It agrees in exact fractions with the direct Gram calculation. The values at $h=1/32,1/256,1/2048$ are $0.41670655,0.40208802,0.40026062$, respectively, rounded. The uniform interval certificate concerns this specified example.

\subsection{Certified finite grids on the half-line}\label{sec:grid}
Quadratic discretization estimates for semi-infinite optimization depend on the form of the constraints and discretization \cite{still}. Stark develops dual bounds to control discretization in inverse problems \cite{stark}. A numerical optimum on a rate grid gives no upper bound for a fraction hidden between its nodes. We give an outward enlargement of the observation box, one-sided in every channel, with a bound independent of the ratio of the latest and earliest sampling times.

In a finite cell $[u,v]$ replace an atom at $k$ by endpoint masses preserving its mean rate. For channel $i$ set $A_i=e^{-t_i u}$, $B_i=e^{-t_i v}$ and $Q_i=(A_i-B_i)/\log(A_i/B_i)$. The maximum increase in its reading is $d_i(u,v)=A_i-Q_i-Q_i\log(A_i/Q_i)$. It occurs where the exponential equals the logarithmic mean $Q_i$, by differentiation. On the terminal cell $[u,\infty]$ move mass to its finite left endpoint, giving $d_i(u,\infty)=A_i$. A convenient bound that avoids cancellation is
\begin{equation}\label{eq:sqgap}
 \tfrac12(\sqrt{A_i}-\sqrt{B_i})^2\le d_i(u,v)\le(\sqrt{A_i}-\sqrt{B_i})^2 .
\end{equation}
The midpoint gives the lower bound. For barycentric weights $p,q=1-p$, the physical exponential at the interpolated point is the weighted geometric mean of $A_i,B_i$. It lies between their weighted harmonic and arithmetic means. The difference of these latter means is $pq(A_i-B_i)^2/(pB_i+qA_i)$. It is bounded by the squared root difference since $pB_i+qA_i-pq(\sqrt{A_i}+\sqrt{B_i})^2=(q\sqrt{A_i}-p\sqrt{B_i})^2$. Reverse Young inequalities give related scalar bounds \cite{young}.

For $N\ge2$ use the union of the finite knots
\begin{equation}\label{eq:knots}
 \left\{-\frac{2}{t_i}\log(j/N):j=1,\ldots,N,\ i=1,\ldots,n\right\}.
\end{equation}
Every channel increases under rounding by at most $N^{-2}$. There are at most $1+n(N-1)$ finite knots. Insert the finite endpoints of a union of closed target intervals, including singletons. Every target atom can now round onto target points. For a general closed target, also insert the minimum and maximum target points of each nonempty finite cell intersection, and the first target point in the terminal cell if present. Round a target atom between these target extrema, or to the first target point of the terminal cell. Other atoms can use the underlying cells. This adds at most two nodes per cell. For infinitely many gaps the input must supply these set extrema. A bounded risk in place of the indicator incurs an additional upper error bounded by its maximum oscillation on any cell, including the terminal cell; this error needs separate control.

In this subsection $A=E$, $P(B,A)=\res_B(E)$, and $\bar P(B,\bar A)=\bar{\res}$ is the compactified value for $T=\bar A$. Write $p_G(B)$ for the mass LP over these finite physical rates with all original bands. Let $d_i$ be an upper bound on the largest channel error over all cells, including the terminal one. One can check the actual nodes instead of requiring exact placement of theoretical knots. Set $B^+=\prod_i[L_i,U_i+d_i]$.
\begin{theorem}[A posteriori enclosure and exact-reading repair]\label{thm:bracket}
Whenever the inner grid LP is feasible,
\begin{equation}\label{eq:grid-sandwich}
 p_G(B)\le P(B,A)\le\bar P(B,\bar A)\le p_G(B^+).
\end{equation}
If $c+b^\top F(g)\ge{\bf1}_A(g)$ at every grid node, an upper bound on the continuous value is
\begin{equation}\label{eq:dualgrid}
 c+\sum_i(b_i^+U_i-b_i^-L_i+b_i^+d_i),\qquad b_i^\pm=\max(\pm b_i,0).
\end{equation}
If $y_0\in B$ and the $n+1$ points $y_0,y_0-R_i e_i$ belong to the grid moment polytope, where $R_i>0$ and $e_i$ are coordinate vectors, set $s=\sum_i d_i/R_i$. Then
\begin{equation}\label{eq:grid-repair}
 p_G(B^+)/(1+s)\le p_G(B)\le p_G(B^+).
\end{equation}
A law giving this lower bound is constructible. Whenever $B\cap\operatorname{int}M\ne\varnothing$, fixed positive $R_i$ work on all sufficiently fine meshes \eqref{eq:knots}, including for an entirely exact moment vector. Thus the value gap is $O(N^{-2})$, with a coefficient depending on the data's position in moment space.
\end{theorem}
\begin{proof}
Rounding a compact optimizer preserves or increases its reward and increases each reading by at most $d_i$. If infinity is rewarded, the physical target is unbounded and has a finite point in the terminal cell. This proves \eqref{eq:grid-sandwich}. Integrating the grid majorant after rounding gives \eqref{eq:dualgrid}.

Take an outer optimizer with readings $y'$. Set $\Delta=(y'-U)_+$ componentwise, so $y=y'-\Delta\in B$. For $s'=\sum_i\Delta_i/R_i>0$, $y_0-\Delta/s'$ is a convex combination of the repair corners. Mix a grid law representing it with the outer optimizer in proportions $s':1$. The readings become $(s'y_0+y)/(1+s')$ and the reward is at least the outer reward divided by $1+s'$. For $s'=0$ nothing needs repair. Since $s'\le s$, this proves \eqref{eq:grid-repair}. In the interior case Theorem~\ref{thm:face} supplies a finite polytope around $y_0$. Approximate its vertices by grid laws: even the first-channel mesh becomes dense on each bounded rate interval. The resulting polytope still contains a smaller fixed neighborhood, including all repair corners.
\end{proof}

The grid size has a sharp worst-case order for this secant construction. A finite cell with error at most $N^{-2}$ spans at most $\sqrt{2}/N$ in each square-root survival coordinate by \eqref{eq:sqgap}; the terminal cell satisfies the same bound. Choose times with successive ratio greater than two. The disjoint physical intervals $[1/t_i,2/t_i]$ each traverse $e^{-1/2}-e^{-1}$ in their own square-root coordinate. Each requires at least $(e^{-1/2}-e^{-1})N/\sqrt2-2$ interior knots. Thus $\Omega(nN)$ knots are necessary in the worst case for this common secant grid.

For fixed times a smallest cell count for prescribed local error tolerances is also obtainable. From each endpoint take the furthest next point satisfying every exact cellwise inequality; stop when the remaining tail can be collapsed within tolerance. Convexity makes the finite-cell errors nondecreasing under interval inclusion. Inductively, each greedy endpoint reaches at least as far as the corresponding point of any admissible mesh, between successive mandatory target endpoints. If the terminal collapse is infeasible, some channel's error limit at a diverging right endpoint exceeds its allowance, so a furthest finite endpoint exists.
\subsection{Exact Hausdorff matrices on a commensurate grid}\label{matrix:hausdorff}
Suppose $t_i=\tau j_i$ with positive integers $j_i$ and $N=\max j_i$. All intermediate powers are retained as free moments. For an interval $D=[a,b]$, $a<b$, write $v_s(x)=(1,x,\ldots,x^s)^\top$ and
\[
 H_{s,g}(y)_{ij}=L_y(x^{i+j}g(x)),\quad 0\le i,j\le s,\qquad L_y(x^j)=y_j.
\]
Let $\mathcal H_N(D)$ consist of all $y=(y_0,\ldots,y_N)$ satisfying positive semidefiniteness of the following blocks:
\[
\begin{array}{c|cc}
N& (s,g)& (s,g)\\ \hline
2d &(d,1)&(d-1,(x-a)(b-x))\\
2d+1 &(d,x-a)&(d,b-x).
\end{array}
\]
Omit a negative-size index. For a singleton use $\mathcal H_N(\{a\})=\{\alpha v_N(a):\alpha\ge0\}$. These are exactly the truncated moment vectors of positive finite measures on $D$. This is the classical Hausdorff criterion (Curto and Fialkow, \href{https://hjm.math.uzh.ch/v017n4/0603CURTO.pdf}{Theorem 4.1 and Remark 4.4}, the latter attributing the even criterion to Kre\u{\i}n and Nudel'man). The zero-mass case follows by continuity, or by separation using the nonnegative-polynomial factorization on a compact interval.

\paragraph{Point formula with singular data.}
For $H\succeq0$ set $d(H,w)=+\infty$ for $w=0$, $d(H,w)=0$ for $w\notin\operatorname{range} H$, and $d(H,w)=1/(w^\top H^\dagger w)$ otherwise. At $z\in I$,
\begin{equation}\label{eq:hausdorff-point}
 A_B(z)=\max_{\substack{y\in\mathcal H_N(I),\,y_0=1\\L_i\le y_{j_i}\le U_i}}\ \
 \min_{(s,g)}d\big(H_{s,g}(y),\sqrt{g(z)}v_s(z)\big).
\end{equation}
Indeed a law with this moment vector can allocate mass $p$ to the anchor precisely when $y-pv_N(z)$ is an interval moment vector (as a submeasure; at maximal $p$ the residual has no anchor atom). Each block becomes $H-pww^\top$. Testing a kernel direction yields the range requirement, and Cauchy--Schwarz on the range gives the stated reciprocal. For the fixed-moment interval extremum, see Schoenberg and Szeg\H{o} \cite{ss} and Schm\"udgen \cite[Theorem 10.29]{schmudgen}. Equation~\eqref{eq:hausdorff-point} recovers the classical maximal point mass for fixed truncated moments, with both support localizers and all intermediate unmeasured powers visible. Equivalently maximize $p$ under the affine block constraints on the residual vector. One global Hankel inverse alone can give a false answer.

\paragraph{Closed unions.}
For $T=\bigcup_{k=1}^{h}D_k\subset I$ a finite union of closed intervals and isolated points, the answer is
\[
 \begin{gathered}\max\ \sum_{k=1}^h y^{(k)}_0,\qquad
 y^{(k)}\in\mathcal H_N(D_k),\quad r\in\mathcal H_N(I),\quad
 \\ y=r+\sum_k y^{(k)},\quad y_0=1,\quad L_i\le y_{j_i}\le U_i.\end{gathered}
\]
Assign the target part of any law to the sets, dividing overlaps by a measurable partition. Conversely each feasible solution yields representing measures on the indicated intervals. Their sum assigns at least the objective to $T$. The compact feasible moment problem attains its maximum. This argument also covers the empty target. It is a direct specialization of classical moment bounds, including \href{https://doi.org/10.1137/S1052623401399903}{Bertsimas and Popescu}. No arbitrary rounding of noncommensurate times to integer powers is justified by this formula.

\subsection{Certificate for the tablet profile}\label{app:tablet}
Write $L_t,U_t$ for the lower and upper bands at time $t$. In $X=e^{-30K}$ write $a=e^{-30(0.008)}$, $b=e^{-30(0.032)}$, $m=L_{30}$ and $v=U_{90}$. A candidate law puts mass $p$ at $a$ and the rest at $d$. Enforcing the two active powers gives
\[
 d=\frac{-a+\sqrt{(a+2m)^2-4(v-m^3)/(a-m)}}{2},\qquad
 p=\frac{m-d}{a-d}.
\]
For $G(x)=(x-d)^2(x+2d)$, the polynomial $G/G(a)$ majorizes the two-tail indicator. Indeed $0<b<d<a<1$, the gap decreases before $d$ and increases after it, and $G(b)/G(a)=1.37272959\ldots>1$. Its linear coefficient is negative and its cubic coefficient positive. The witness saturates the correct box sides and has second and sixth powers $0.38411869\ldots$ and $0.05981127\ldots$, both strictly inside their bands. Thus the globally sharp maximum is
\[
 \boxed{\quad p=0.039701199121\ldots = 3.9701199121\ldots\%\quad}
\]
and the unique extremizer allocates its entire outlying fraction to the slow side. Strict inequality on the fast side rules out positive fast mass in an optimizer. This finding concerns the distribution that maximizes the specified joint tail. Other feasible fits may contain fast material. The allowance is a deterministic modeling choice; it does not have a replicate-level coverage interpretation. Independent finite-grid verification at $N=256$ encloses the maximum between $3.9647\%$ and $3.9860\%$ (rounded outward). A separate arithmetic certificate encloses the exact formula above.

\subsection{Computing a certificate for a new experiment}\label{app:practical}
For an interval or a union of intervals, insert the finite target endpoints into a common rate grid for all channels. Let $p_G(B)$ be the largest grid fraction under the original bands. If a grid population is feasible, it supplies an attainable lower bound. For an upper bound replace any rate between two adjacent finite grid nodes by an interpolation at the endpoints preserving its mean rate. Convexity of the exponential increases each reading. The largest increase is at most the squared difference between the square roots of the endpoint signals. A terminal cell controls arbitrarily fast rates. Target endpoints allow this rounding to preserve or increase the fraction in the target.

Taking knots $-2\log(j/N)/t_i$ for $1\le j\le N$ in each channel puts every reading increase below or at $N^{-2}$. Thus an outer grid calculation permits the enlarged upper endpoints $U_i+N^{-2}$, while retaining the measured lower endpoints. A separate verifier checks the achieved distribution by interval arithmetic and constructs a continuous upper bound from the signed coefficients of the grid majorant. If the observation box meets the interior of the moment body, a correction among finitely many grid distributions restores the original bands at a loss of order $N^{-2}$ in target mass, including for exact readings. Constants depend on the position of the data in the moment body. There are at most $1+n(N-1)$ basic finite nodes before target endpoints are added. For $q$ nodes the verifier builds $nq$ channel enclosures. Exact weight reconstruction also requires solving a subsystem on the proposed support; wide intervals near a singular basis may prevent certification. Supplement, Section~\ref{sec:grid} proves the bracket and gives its data-dependent constant.

The method can return a conditional upper bound when a physical witness has not been certified. For instance, $\mathbb E e^{-K}=1/2$ and $\mathbb E e^{-2K}=1/4$ force $K=\log2$ by zero variance; a grid of rational physical rates cannot reproduce this exact pair. The implementation handles this family by a separate symbolic certificate. It accepts a numerical optimizer's output only after checking the continuous upper bound and any claimed finite lower witness.
\paragraph{Three irregularly spaced readings.}
For the numerical verifier take $t=(1/2,\sqrt2,3)$ and the target rate intervals $[0.05,0.11]\cup[1.8,2.05]$. We tested both the box $[0.540,0.546]\times[0.262,0.269]\times[0.138,0.144]$ and the exact moment vector $(0.543,0.265,0.141)$. Table~\ref{tab:brackets} gives enclosed continuous maxima on the physical half-line. The inner witness specifies finite physical nodes and an exactly imposed linear subsystem for its weights; interval Gaussian elimination verifies positivity and the remaining bands. The outer certificate uses signed coefficients, their grid values and the one-sided rounding errors. Numerical optimization only proposes the witnesses.
\begin{table}[htbp]\centering
\begin{tabular}{@{}rrr@{}}\toprule
Mesh parameter $N$ & Three two-sided bands & Three exact readings\\\midrule
16  & $[0.80603,0.84465]$ & $[0.73415,0.78500]$\\
32  & $[0.80623,0.81888]$ & $[0.73470,0.74736]$\\
64  & $[0.80628,0.80945]$ & $[0.73484,0.73801]$\\
128 & $[0.80632,0.80712]$ & $[0.73491,0.73571]$\\\bottomrule
\end{tabular}
\caption{Outward-rounded enclosures on the total population fraction. Each endpoint is checked with exact rational interval operations and enclosed exponentials.}\label{tab:brackets}
\end{table}

\section{Why a lifetime window can hold more than any point in it}
A detector with a fixed lifetime bandwidth counts emission from every lifetime in its window. Maximizing over the possible position of one hidden lifetime answers a different question. The distinction persists with four exact readings and a bounded physical rate range.

\paragraph{A four-reading counterexample.}
Put $X=e^{-\tau K}=1/2+Y/5$ and restrict $Y$ to $[-2,2]$. Measurements at $\tau,2\tau,3\tau,4\tau$ determine the first four moments of $Y$ by an invertible affine change of basis. Set them to $0,1,0,31/16$. For the closed window $J=[-1/2,1/2]$ in $Y$, the largest compatible fraction is exactly $5/8$ and the smallest is zero. The largest mass compatible with any single point in $J$ is exactly $20/37$. The interval and single-point maxima differ by $25/296$.

The polynomial $Q(y)=(y^2-9/4)^2/4$ is nonnegative everywhere and at least one on $J$. Its expectation is $5/8$. Put total weight $5/8$ at $\pm1/2$ and the remaining $3/8$ at $\pm3/2$, splitting each pair equally. This reproduces the moments and attains the upper bound. To attain zero, place total weight $240/289$ at $\pm3/4$ and $49/289$ at $\pm\sqrt{22/7}$. Both pairs are outside $J$ and give the same moments.

For a single point, $1,Y,(4/\sqrt{15})(Y^2-1)$ are orthonormal under these moment data. Multiply each polynomial by its value at the query $y$, sum, and divide by the sum of the three squared evaluation values. The resulting polynomial equals one at $y$. Its mean square is $[1+y^2+(16/15)(y^2-1)^2]^{-1}$ and bounds the mass at $y$. Writing $u=y^2$, the excess of the denominator over $37/20$ factors as $(1/4-u)(13/15-16u/15)$, which is nonnegative for $0\le u\le1/4$. All witnesses correspond to finite positive physical rates.

\paragraph{The three atoms attaining the point bound.}
Let $r=\sqrt{601}$. Assign mass $20/37$ to the query $y=1/2$. Place the remaining mass at
\[
 y_\pm=(5\pm r)/16,\qquad w_\pm=17/74\mp245/(74r).
\]
Since $24<r<25$, the two sites satisfy $-2<y_-<-1/2<1/2<y_+<2$.
Both weights are positive since $17r>245$. The three sites therefore give $0<X<1$ and finite positive physical rates.

To check the moments, let $b_j=w_-y_-^j+w_+y_+^j$. Direct addition gives
$b_0=17/37$ and $b_1=-10/37$. Each outer site satisfies $y^2=(5/8)y+9/4$.
Consequently $b_{j+2}=(5/8)b_{j+1}+(9/4)b_j$ and
$(b_2,b_3,b_4)=(32/37,-5/74,1127/592)$.
Adding the queried atom contributes $(20/37)2^{-j}$ to each $b_j$.
The resulting moments for $j=0,1,2,3,4$ are $1,0,1,0,31/16$, including normalization.
This attains the pointwise upper bound at the right endpoint of $J$. Reflection gives the same mass at its left endpoint. The factored denominator bound is strict in the interior.
It follows that $5/8-20/37=25/296$ of the total signal can be assigned to the window
over and above what any single point there can hold.

\paragraph{When can a point certificate be widened?}
Suppose a verified nonnegative bounding function $q$ has expectation at most $U$ throughout the reading box. If $q\ge m>0$ on a closed interval $J$, then every compatible distribution satisfies $\mu(J)\le U/m$. Integrate $m{\bf1}_J\le q$. If a physical witness puts mass $L$ at a point in $J$, the interval maximum lies in $[L,\min(1,U/m)]$. For polynomial $q$ on consecutive grids, root tests or Bernstein subdivision give a certified $m$. In particular, if $U/\omega$ and $L/\omega$ both converge to the same finite value and $m\to1$, then the normalized interval maximum has that limit as well. Small gaps between reference sites can make $q$ vary rapidly across the window. Its verified minimum supplies a configuration-specific test.

A guaranteed minimum in an open window instead comes from maximizing mass in its closed complement. Classical interval moment theory and the equal-weight two-rate central-gap limit explain the ingredients of this special case. A sharp law for intervals crossing moving, unequally weighted reference sites at arbitrary scales would require control of multiple contact points. The accompanying exact-arithmetic and Lean files identify the computations and certificate implications checked by machine.

\section{Can a decay reveal heterogeneity without locating it?}
Consider a known two-rate population with one-step survivals $1/2-g$ and $1/2+g$ and fast amplitude $w$. No component lies in the test window $0.45\le X\le0.55$ for the separations in Figure~\ref{fig:phase}. We compute the largest mass that \emph{any} compatible distribution can put in that window. For $w=1/4$ and a reading half-width of $0.2\%$, increasing the survival separation from $0.20$ to $0.60$ reduces this maximum from 98.36\% to 41.68\%. At both separations, every single-rate signal can already be distinguished from the two-rate reference at the stated error. Thus a measurement can establish heterogeneity while leaving most of the signal free to move into a window containing neither true rate.

\begin{figure}[htbp]\centering
\includegraphics[width=\linewidth]{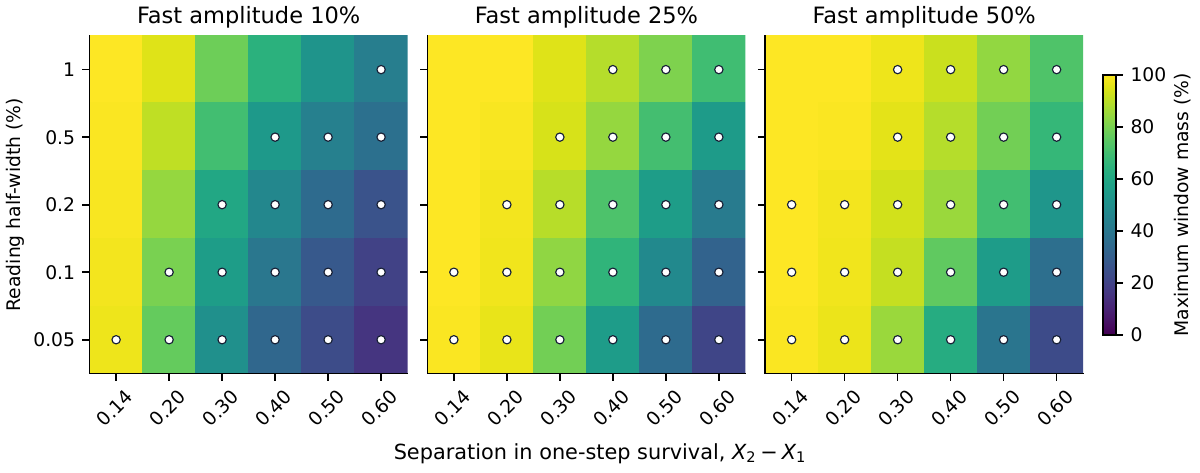}
\caption{Finite-precision resolution for four readings at times $1,2,3,4$. Color gives the maximum admissible amplitude in the central window. The certified upper and lower values differ by less than 0.005 percentage points in every cell. A dot means the error boxes of the reference and every possible single-rate curve are disjoint, both at the indicated half-width. The gap between two survival coordinates is dimensionless and converts to rates by $K=-\log(X)/\Delta t$.}\label{fig:phase}
\end{figure}

A separate fit comparison uses $X=0.3,0.7$, $w=1/4$, half-width $0.002$, and a fixed perturbation $(0.35,-0.50,0.25,-0.20)$ times that half-width on the four readings. NNLS and a free two-rate fit assign zero to the window. A first-difference Tikhonov fit assigns 18.98\%. All three fits satisfy the same error tolerance. The certified maximum lies between 70.0924\% and 70.0967\% (outward rounded). Convex mixtures of the physical witness and the true population realize all intermediate fractions up to the witness value. On grids of 100 to 400 points the smoothed fit varies from 19.76\% to 18.61\%. The bound includes explicit checks between grid nodes and on the whole compactified axis. In the optical case all three fit types instead agree to within a small fraction of a percentage point on the higher-count draw, assigning roughly 98\% to the $1$--$8$ ns window. The data alone guarantee about 54.49\%. Agreement among algorithms and resolution of the underlying population are distinct tests.

For the same two rates and exact reference centers, one further reading can be planned before perturbing the data. Give every candidate time the same absolute half-width $0.002$. Among the 16 additional integer times 5 through 20, time 11 minimizes the maximum compatible central-window fraction. It reduces the four-reading value of about 72.14\% to about 43.89\%. Time 10 permits more than 43.92\%, and time 12 more than 44.62\%. Upper and lower certificates separate time 11 from each competitor. This finite-menu result assumes reference-centered future readings and constant achievable absolute precision. It complements the small-gap timing and tenth-power examples in the technical material.

\section{Additional limiting regimes and certified examples}
These derivations accompany the article. In \eqref{eq:supp-nested-summary}, $A_h$ denotes the nested three-rate example. In general $A$ is the maximum allowed mass at the queried rate, $S=\omega^{-1}$, $d^2=P(z)^2/S$, $\eta$ is the raw half-width multiplier, and $C$ denotes the design coefficient in the single-cluster theorem. The finite envelope below pertains to the nested three-rate data at $\eta=h^{10}/728$.
\begin{equation}\label{eq:supp-nested-summary}
 2/5<A_h<2/5+h,\qquad 0<h\le1/32 .
\end{equation}
\paragraph{A critical exponent of ten from three rates.}
For times $1,\ldots,6$ with equal widths, take $k_*=\log2$, equal weights at displacements $0,2\log(1+h^2),\log(1+h)$, and anchor $z=\log(1+h^2)$. Here $S\to3/2$, $D^2/h^{10}\to1$, and $C=728$. At $\eta/h^{10}\to\lambda>0$ the sharp limit is $(3/2+1/(728\lambda))^{-1}$. At $\eta\ll h^{10}$ the limit is zero by monotonicity. Throughout $h^{10}\ll\eta\ll h^2$ the limit is $2/3$. For the lower bound send $\lambda$ to infinity at the finer transition. For the upper bound, with $X=e^{-K}$ let $a=e^{-(k_*+\log(1+h))}$ and $b=e^{-(k_*+z)}$. The globally nonnegative quadratic $(X-a)^2/(b-a)^2$ is one at the anchor, costs $2/3+O(h)$ on the reference and amplifies the reading errors by $O(h^{-2})$. Thus even error much larger than $h^{10}$ can leave a macroscopic portion of the total mass unavailable at the anchor. When $\eta\gg h^2$ the limit is one: an interior rate at distance $O(h)$ from the anchor and mass $O(1)$ can match its first displacement moment. More explicitly, for fixed $\epsilon\in(0,1)$ give $1-\epsilon$ to the anchor, and $\epsilon$ to the unique rate matching the mean displacement of the reference. Its displacement is $O(h/\epsilon)$ and it is physical for small $h$. Taylor's theorem bounds every resulting reading discrepancy by $O_\epsilon(h^2)$. Send $\epsilon\downarrow0$.

At the intermediate coarse transition $\eta/h^2\to\lambda>0$, the scaled reference tends to mass $2/3$ at zero and $1/3$ at one, and the scaled anchor tends to zero. With $C_{*,2}=4$ from the same six readings, the limiting anchor mass is
\[
 \left(1+\frac{1/9}{2/9+4\lambda}\right)^{-1}.
\]
Indeed, use the degree-two construction with the ordinary scaled first two moments of this reference. Their limiting values both equal $1/3$; adding the sharp positive $4\lambda$ to the second gives the stated degree-one Christoffel kernel at zero. Its monic orthogonal polynomial $x-1/3$ has nonzero value at zero. The degree-one anchored quadrature has positive weights at the anchor and one simple distinct finite node, so its Jacobian stays nonsingular under the moving moments and anchor. Taylor recovery at all six times and a two-time global double-zero majorant give the same upper and lower limits as in the common-scale polynomial reduction. Thus the transition powers in this example are ten and two, with no transition at six or four.

\paragraph{Two separated tails on an irrational time grid.}
Take $t_1=1$, $t_2=\rho=\sqrt 2$ and $A=[0,a]\cup[b,1]$, $a=1/4$, $b=3/4$.
The observed bands are $[0.50,0.51]$ and $[0.38,0.39]$.
Put $m=(b^\rho-a^\rho)/(b-a)$ and $c=(m/\rho)^{1/(\rho-1)}$.
Strict convexity gives $a<c<b$. The tangent gap
$g(x)=x^\rho-c^\rho-m(x-c)$ is minimized at $c$, and $g(a)=g(b)>0$.
Thus $Q=g/g(a)$ majorizes the two-tail indicator with exactly three contacts $a,c,b$.
The maximum fraction and realizing weights are
\[
 p=\frac{0.39-c^\rho-m(0.50-c)}{g(a)},\qquad
 w_c=1-p,\quad w_b=\frac{0.50-c-p(a-c)}{b-a},\quad w_a=p-w_b .
\]
All three weights are positive and give moments $0.50$ and $0.39$ by substitution.
The two nonzero coefficients and affine independence give uniqueness.
Numerical enclosures and a separate original-band linear program are supplied in the archive.
All three optimizing physical rates are finite.

The directed calculation gives $p=0.52829535744388\ldots$; a mesh with $12\,801$
points that omits $c$ gives $0.528295356697$. As a sharp one-interval example take physical
survival support $[1/5,1]$, times $1,2,3,4$, and target $[2/5,4/5]$. Write
$Y=(5x-3)/2$ and $Q=1+(1/4-Y^2)Y^2/(3/4)$. It majorizes the target indicator,
with contacts exactly at $x=i/5$, $i=1,\ldots,5$. Take the equal-weight moments there,
and place each at its corresponding band endpoint according to the sign of the coefficient
of $Q$, with width $1/1000$ on the other side. The sharp fraction is $3/5$ and the
Vandermonde determinant is $288/9765625$. Hence every optimizing law uses all five
rates with equal weights. All moment bounds and this certificate are checked in rational
arithmetic in the supplied program.
\paragraph{Checks in exact arithmetic.}
For $I=[0,1]$, exact first and second moments $1/2,2/5$, and $z=1/2$, the ordinary Gram bound equals one. The interval localizer gives $(1/2-2/5)/(1/4)=2/5$. It is attained by masses $3/10,2/5,3/10$ at $0,1/2,1$. In physical finite rates this supremum is approached but cannot be attained. A disconnected-target example has bands $1/2\le y_1\le51/100$, $27/100\le y_2\le9/32$ and target $[0,1/4]\cup[3/4,1]$. The sharp mass is $1/2$, by the certificate $16(x-1/2)^2$ and masses $1/4,1/2,1/4$ at $1/4,1/2,3/4$. The supplied generic solver uses exact rational-polynomial feasibility with singular semidefinite blocks allowed. It returns certified dyadic brackets and, on request, recovers an algebraic atomic lower witness independently from the block computation. Finite time limits can leave a decision unresolved. Principal-minor encoding grows exponentially in block dimension, so the decisive parameter for this implementation is the largest integer power after common-unit reduction.

A sparse-time example on $[1/5,1]$ uses powers $(1,2,4,5)$, anchor $2/5$,
and masses $(1,3,1,2)/7$ at $(1,2,3,5)/5$. Take these readings as the upper,
lower, upper and lower bounds, respectively, widening each band by $1/100$ on the
other side. The certificate $(6250/37)(x-1/5)(1-x)(x-3/5)^2(x+5/6)$ has
zero cubic coefficient, equals one at the anchor, and gives the exact value $3/7$.
All four nodes are necessary.
\section{Alternating certificates for a prescribed rate}\label{p4:alternating}
Put $x=e^{-\tau K}$ for any time unit $\tau>0$. Write $I=[l,u]\subset[0,1]$, $l<u$, with $l=0$ allowing the infinitely fast endpoint. Let $0=\lambda_0<\lambda_1<\cdots<\lambda_n$, with $\lambda_i=t_i/\tau$ and $x^0=1$ everywhere. Let the finite closed box $B=\prod_{i=1}^n[L_i,U_i]$ meet the normalized moment body $M=\operatorname{conv}\{(x^{\lambda_1},\ldots,x^{\lambda_n}):x\in I\}$. For $Q=\sum c_i x^{\lambda_i}$ write
\begin{equation}\label{eq:supp-box}
 F_B(Q)=c_0+\sum_{c_i>0,\ i>0}c_iU_i+\sum_{c_i<0,\ i>0}c_iL_i.
\end{equation}
At $z\in I$ seek $A_B(z)=\max_{\mu:\,m(\mu)\in B}\mu(\{z\})$, with probabilities on $I$. Compactness and upper semicontinuity give attainment.

The alternating geometry of moment boxes is classical. Kre\u{\i}n and Nudel'man, \href{https://doi.org/10.1090/mmono/050}{Chapter VI, Theorems 1.1, 4.4 and 5.1--5.3}, analyze alternating supporting normals, motion of prescribed maximal atoms, and extrema for one-sided integral objectives. We give a self-contained selected-row formula for power readings, including singular boxes.

\begin{theorem}[Alternating subgrid formula]\label{matrix:alternating}
For every such box, including zero widths, and any prescribed point including an endpoint, the maximizing law is unique. The sharp value is the following finite-dimensional minimum.

Select a nonempty subset $J\subset\{0,\ldots,n\}$, $m=|J|$. Choose a multiset of $m-1$ roots in $I\cap(0,\infty)$ away from $z$; interior roots have even multiplicity. In the selected span require all these zeros with their multiplicities. The nullspace is one-dimensional. Take the nonnegative orientation and normalize at $z$ whenever its value there is positive. At $l=0$, there is also a compulsory zero there when $0\notin J$. Then
\begin{equation}\label{eq:supp-subgrid}
 A_B(z)=\min_{J,\ {\rm roots}} F_B(Q),\qquad Q(z)=1 .
\end{equation}
The nonzero coefficients strictly alternate in increasing exponent order. By \eqref{eq:supp-box}, only alternating band endpoints on the selected measurement rows enter a certificate. The optimizer is supported on $z$ and the zeros of an optimal certificate. Its weights solve linear equations on these nodes using normalization and the selected endpoint readings. For $n\ge1$ at most $\lfloor n/2\rfloor+2$ rates suffice. This uniform count is sharp on compact physical rate intervals.
\end{theorem}
\begin{proof}
To prove \eqref{eq:supp-subgrid}, first observe that $A_B(z)$ is the infimum of $F_B(Q)$ over $Q\ge0$, $Q(z)=1$. The convex set $\{(m(\mu),s):s\le\mu(\{z\})\}$ is closed by weak compactness. Separate it strictly from the compact $B\times\{A_B(z)+\epsilon\}$. The coefficient on $s$ is nonnegative by downward closure and is positive since the moment projections intersect. Division gives an affine majorant of the indicator (test point masses) with box upper expectation below $A_B(z)+\epsilon$. Scaling down to $Q(z)=1$ preserves the bound since $F_B(Q)\ge0$ for nonnegative $Q$ by feasibility. The converse follows by integration.

We show dual attainment without requiring positive widths. If $B$ contains a point $y$ interior to $M$, its expectation functional $L_y$ is strictly positive on every nonzero nonnegative $Q$. A nonconstant $Q$ vanishing in expectation there would give a supporting hyperplane through an interior point. Compactness on a coefficient unit sphere gives $F_B(Q)\ge L_y(Q)\ge c\|\operatorname{coeff}(Q)\|$ on the nonnegative cone for some $c>0$. This bounds a minimizing sequence. The body has full dimension by linear independence.

Otherwise separate $B$ from the open interior of $M$. Since they meet in closure, there is a nonconstant $Q_0\ge0$ with $F_B(Q_0)=0$. Every feasible law is supported on its zeros, and every nonzero measured coefficient fixes its individual box endpoint. Let $J$ be its nonzero coefficient support, $\lambda_*=\min_{i\in J}\lambda_i$, and write $Q_0=x^{\lambda_*}R$. The positive zeros have total multiplicity $k\le |J|-1$ by the exponential zero bound (in logarithmic coordinates, divide by the first exponential, differentiate, and apply Rolle with multiplicities). At $l=0$, $R(0)\ne0$. For $z$ off the zero set, $Q_0/Q_0(z)$ gives value zero. For $z$ on it, interpolate $V$ on the positive zeros with values $\mathbf1_{\{z\}}$ and with derivatives zero up to one less than each multiplicity. The $k$ Hermite conditions on the selected powers have full row rank since any $k$ distinct selected powers have nonsingular interpolation matrix by the same zero bound. If $l=0$ and $\lambda_*>0$, include the fixed independent constant $\mathbf1_{\{z=0\}}$ and solve the same conditions in the selected powers. Otherwise take $V$ entirely in the selected span. Near a positive nonanchor zero, $V/Q_0$ stays bounded. Near a zero anchor $V$ is positive. At a compulsory zero at the origin, either $V(0)=1$ or $V=O(x^{\lambda_*})$. Hence $V+D Q_0\ge0$ on $I$ for all large $D$. Its measured coefficients retain the signs of $Q_0$. Every feasible measure integrates it to its box upper value, which equals the mass at $z$. This gives attainment also at a boundary box. In this case there is only one feasible law: its weights on the positive zeros are determined by selected readings and normalization (adjoin the constant row for a possible atom at zero). Generalized Vandermonde matrices at distinct positive nodes have full column rank when they have enough rows by the zero bound. A node at zero is separated by the constant row and the positive powers.

Intersect the nonnegative cone with the attained certificate's closed coefficient orthant, also splitting the constant sign. There $F_B$ is linear. A compact base obtained by normalizing Lebesgue integration on $I$ shows every element is a finite sum of extreme rays. Each ray $P$ has $F_B(P)\ge0$, and if $P(z)>0$, $F_B(P)/P(z)\ge A_B(z)$. At least one ray with positive evaluation thus attains the optimal ratio.

In an extreme ray's nonzero support $J$ factor out the smallest power. If the remaining positive zeros on $I$ have total multiplicity below $m-1$, there is an independent polynomial satisfying their homogeneous Hermite conditions in the same span. Its quotient by the ray stays bounded on $I$ (factor out the same power at zero). Adding and subtracting small multiples preserves nonnegativity and coefficient signs, contradicting extremality. Thus there are exactly $m-1$ positive zeros counted with multiplicity. The interpolation kernel has dimension one. Generalized Descartes' rule implies at least $m-1$ changes of coefficient sign. For completeness, divide by the lowest power on $x>0$ and differentiate; the derivative has the signs of the remaining coefficients. If the first two signs agree and a positive zero exists, the derivative has a further zero before the first one, because the function starts by moving away from zero. Rolle's theorem at the remaining zeros with multiplicities completes induction. Opposite first signs allow the single extra zero.

Complementary slackness confines all maximizing laws to the zeros and $z$ and fixes every selected reading at its indicated endpoint. For $l>0$ there are at most $m-1$ distinct zeros. For $l=0$ with smallest exponent positive there may be one additional zero; adjoining normalization gives one additional row. Generalized Vandermonde rank gives unique weights, also when $Q$ is constant. Counting interior zeros twice and positive endpoints once yields the sharper support bound. For $l>0$ there are at most $\lfloor(m+1)/2\rfloor$ distinct zeros. For $l=0$ there are at most $\lfloor m/2\rfloor$ positive zeros, and one extra zero only when $m\le n$ and the constant is absent. Add the prescribed point in each case.

For sharpness take $l>0$ and all $n+1$ powers. For $n=2d$ place simple zeros at both endpoints and double zeros at $d-1$ interior points. For $n=2d+1$ take one simple endpoint zero and $d$ double interior zeros. The resulting interpolant has fixed sign throughout $I$ and no other zero. Add a distinct target $z$. Give positive weights to all these nodes. Set each resulting reading at the band endpoint selected by the certificate and widen the band on the opposite side. These data attain the bound, and uniqueness forces every one of the $\lfloor n/2\rfloor+2$ atoms.
\end{proof}
The certificate is a determinant with root evaluations and derivatives as rows and the evaluation at $x$ as last row, divided by its value at $z$. The number of continuously variable interior root locations is at most $\lfloor n/2\rfloor$. The number of candidate subgrids can be exponential. Unselected readings can place the optimum strictly inside some bands. For example $\delta_z$ may fit strictly inside the whole box.

For a finite anchor on the physical half-open interval, Theorem~\ref{thm:face} identifies this compact value with the physical supremum whenever the physical problem is feasible. By uniqueness, attainment is equivalent to zero mass at the infinitely fast endpoint in the compact optimizer.

\section{Checking a product-yield observable}
For PET photoreforming, Qin et al.\ report soluble terephthalic-acid-equivalent conversion at eight times \href{https://doi.org/10.1021/jacs.5c08549.s001}{(supporting data, Table S7)}. If one interprets its complement as a normalized static-rate signal, the reported means at 8 and 20 hours give $s(8)=0.50$ and $s(20)=0.01$. Jensen's inequality requires $s(20)\ge s(8)^{5/2}$. The smallest uniform allowance for all eight means is $9.4626$ percentage points, with the two indicated times binding. In the four-hour survival coordinate, the unique root $u\in(0,1)$ of $u^2+u^5=0.51$ gives a single rate attaining that discrepancy; the other six readings have strict slack (supplementary certificate). The measured product yield and remaining polymer are distinct observables. This example tests the proposed observation model for product formation and by itself assigns neither rate fractions nor a reaction mechanism to the experiment.

\section{Transition-law benchmarks}
\paragraph{Three nearby rates and one extra reading.}
Take six equally precise readings at times $1,\ldots,6$ and $k_*=\log 2$ in reciprocal time units. Give equal weights to three physical displacements $0$, $2\log(1+h^2)$ and $\log(1+h)$, and query the intermediate displacement $\log(1+h^2)$. One reference gap is of order $h^2$, while the outer rate is of order $h$ away. The anchor has zero weight in the reference. Nevertheless $S\to3/2$, $D^2/h^{10}\to1$, and the grid-design formula of the main article gives $C=728$. At tolerance $\eta=h^{10}/728$ the largest compatible anchor fraction tends to $2/5$.

The finite-$h$ statement is stronger than a fitted slope. For every $0<h\le1/32$ an integer-polynomial certificate gives $2/5<A_h<2/5+h$. At $h=1/32$ the maximum is $0.4167065509$ to ten decimal places. The raw moment matrix at this scale has infinity-norm condition number $6.79\times10^{17}$. Rounding its exact entries to binary64 and solving directly gives $0.2587118847$ on the recorded implementation. Assembling the normalized form from the reference parameters and its structured moment update gives $0.4167065509$, with condition number about $1.68$. These are different input scenarios. Re-expressing already rounded raw moments cannot recover discarded information. The exact checks in the finite-scale certificate in Section~\ref{app:nested-examples} use rational arithmetic.

There is a concrete design consequence. Keep the same tolerance at every measurement and add a reading at time seven. Minimizing the design cost over the seven possible six-row subsets now gives $C=2132/7$ at $(1,2,3,4,6,7)$. The new limiting hidden fraction is
\[
 \left(\frac32+\frac{728}{2132/7}\right)^{-1}=\frac{82}{319}=0.25705\ldots .
\]
The seventh reading excludes about fourteen additional percentage points of hidden mass in this asymptotic experiment. The reference rates and the precision of each instrument reading have stayed fixed at each $h$. No claim of global optimality over all additional times is needed for this comparison. In general, a candidate time can be evaluated by recomputing this design cost, with its achievable precision reflected in $\sigma_j$.

A weak component changes both the tolerance and the mass scale. With two reference rates at displacements $0,2\log(1+h)$, weights $h,1-h$, four consecutive readings with unit widths, and query $\log(1+h)$, one has $\omega=4h(1-h)$, $d^2\sim4h^5$ and $C=80$. At tolerance $h^5/20$, the largest hidden fraction is $2h+o(h)$. In the three-rate example it stays of order one.

\paragraph{Timing a measurement of a suspected split.}
If the competing physical models are specified in advance, their useful measurement times can sometimes be found without solving an extremal problem. Compare a rate $k$ with an equal split at $k-s,k+s$, where $0<s<k$. The difference in surviving signal is
\[
 \Delta(t)=e^{-kt}\{\cosh(st)-1\}.
\]
With the same absolute precision available at all times, its maximum occurs at
\[
 t_{\rm pair}=\frac{2}{s}\operatorname{atanh}(s/k)\ \longrightarrow\ 2/k
 \quad (s/k\to0,\ k\ \text{fixed}).
\]
Indeed, the derivative at $t>0$ has the sign of $s-k\tanh(st/2)$. Thus a reading near two central lifetimes is most informative for this particular pair when the split is small. To exclude one prediction for every allowed reading from the other, the signal difference must exceed twice the absolute half-width. The unrestricted population question uses the resolution function of the main article.

\paragraph{Two separated pairs.}
For the separated-rate theorem of the main article, take survival coordinates
$x=(1/4,1/4+h,3/4,3/4+h)$, $z=1/4+h/2$ and weights $(1,2,3,4)/10$.
Use $\eta/d^2=1/16$ and $\sigma=(1,2,1,2,1,2,1,2)$. The product formula gives $\epsilon=O(h)$ even though the two centers remain half a survival unit apart. The entries in Table~\ref{tab:separated} are rounded from exact fractions. For $h=1/128$ the absolute normalized error decreases from $9.04\times10^{-4}$ to $4.32\times10^{-6}$ when the signed correction is included.
\begin{table}[htbp]\centering
\begin{tabular}{@{}lrrr@{}}\toprule
$h$ & Exact $A/\omega$ & $f(\xi)$ & Including correction\\\midrule
$1/64$ & 0.689686488 & 0.691356788 & 0.689704726\\
$1/128$ & 0.685900931 & 0.686804748 & 0.685905249\\\bottomrule
\end{tabular}
\caption{Two separated pairs of rates. The exact values use both the admissible moment vertex and a positive physical quadrature.}\label{tab:separated}
\end{table}

\endgroup

\begin{thebibliography}{99}
\setlength{\itemsep}{0.25em}

\bibitem{prony}
A. Akinshin, G. Goldman and Y. Yomdin (2021).
\href{https://doi.org/10.1090/mcom/3571}{Geometry of error amplification in solving the Prony system with near-colliding nodes}.
\emph{Mathematics of Computation} \textbf{90}(327), 267--302.

\bibitem{fourier}
D. Batenkov, G. Goldman and Y. Yomdin (2021).
\href{https://doi.org/10.1093/imaiai/iaaa005}{Super-resolution of near-colliding point sources}.
\emph{Information and Inference: A Journal of the IMA} \textbf{10}(2), 515--572.

\bibitem{bp}
D. Bertsimas and I. Popescu (2005).
\href{https://doi.org/10.1137/S1052623401399903}{Optimal inequalities in probability theory: a convex optimization approach}.
\emph{SIAM Journal on Optimization} \textbf{15}(3), 780--804.

\bibitem{cf}
R.~E. Curto and L.~A. Fialkow (1991).
\href{https://hjm.math.uzh.ch/v017n4/0603CURTO.pdf}{Recursiveness, positivity, and truncated moment problems}.
\emph{Houston Journal of Mathematics} \textbf{17}(4), 603--635.

\bibitem{dip2}
P. Dahinden (2021).
\href{https://github.com/cran/disprofas/blob/c7aebfa311fb5d29cab7ff09a85f0ce2f5f431d0/man/dip2.Rd}{Dissolution data of one reference batch and five test batches}.
Dataset \texttt{dip2} in version 0.1.2 of the R package \href{https://doi.org/10.32614/CRAN.package.disprofas}{\emph{disprofas: Non-Parametric Dissolution Profile Analysis}}, Comprehensive R Archive Network. Source data from Shah et al., Table~4.

\bibitem{didio}
P.~J. di Dio and K. Schm\"udgen (2022).
\href{https://doi.org/10.1016/j.jmaa.2022.126066}{The multidimensional truncated moment problem: the moment cone}.
\emph{Journal of Mathematical Analysis and Applications} \textbf{511}(1), article 126066.

\bibitem{young}
S. Furuichi and N. Minculete (2021).
\href{https://doi.org/10.3390/e23050514}{Refined Young inequality and its application to divergences}.
\emph{Entropy} \textbf{23}(5), article 514.

\bibitem{vander}
W. Gautschi (1962).
\href{https://doi.org/10.1007/BF01386302}{On inverses of Vandermonde and confluent Vandermonde matrices}.
\emph{Numerische Mathematik} \textbf{4}, 117--123.

\bibitem{jameson}
G.~J.~O. Jameson (2006).
\href{https://doi.org/10.1017/S0025557200179628}{Counting zeros of generalised polynomials: Descartes' rule of signs and Laguerre's extensions}.
\emph{The Mathematical Gazette} \textbf{90}(518), 223--234.

\bibitem{ks}
S. Karlin and W.~J. Studden (1966).
\emph{Tchebycheff Systems: With Applications in Analysis and Statistics}.
Pure and Applied Mathematics, vol.~15. Interscience Publishers, New York.

\bibitem{krein}
M.~G. Kre{\u\i}n and A.~A. Nudel'man (1977).
\href{https://doi.org/10.1090/mmono/050}{\emph{The Markov Moment Problem and Extremal Problems: Ideas and Problems of P.~L. \v{C}eby\v{s}ev and A.~A. Markov and Their Further Development}}.
Translations of Mathematical Monographs, vol.~50. American Mathematical Society, Providence, RI.

\bibitem{le}
D. Le, T. Pham, A. Rinaldo and N. Ho (2026).
\href{https://doi.org/10.48550/arXiv.2610.25857}{Weak limits of parameter estimation in strongly identifiable mixture models}.
Preprint, \href{https://arxiv.org/abs/2610.25857v1}{arXiv:2610.25857v1}.

\bibitem{liu}
P. Liu and H. Ammari (2023).
\href{https://doi.org/10.1093/imaiai/iaad048}{Super-resolution of positive near-colliding point sources}.
\emph{Information and Inference: A Journal of the IMA} \textbf{12}(4), 3087--3111.

\bibitem{shah}
V.~P. Shah, Y. Tsong, P. Sathe and J.-P. Liu (1998).
\href{https://doi.org/10.1023/A:1011976615750}{In vitro dissolution profile comparison---statistics and analysis of the similarity factor, $f_2$}.
\emph{Pharmaceutical Research} \textbf{15}(6), 889--896.

\bibitem{simon}
B. Simon (2008).
\href{https://doi.org/10.1090/pspum/079/2500498}{The Christoffel--Darboux kernel}.
In D. Mitrea and M. Mitrea (eds.), \emph{Perspectives in Partial Differential Equations, Harmonic Analysis and Applications}, Proceedings of Symposia in Pure Mathematics \textbf{79}, American Mathematical Society, 295--335.

\bibitem{still}
G. Still (2001).
\href{https://doi.org/10.1007/s101070100239}{Discretization in semi-infinite programming: the rate of convergence}.
\emph{Mathematical Programming} \textbf{91}(1), 53--69.

\bibitem{watson}
G. A. Watson (1990).
\href{https://doi.org/10.1093/imanum/10.4.569}{Chebyshev approximation to data by positive sums of exponentials}.
\emph{IMA Journal of Numerical Analysis} \textbf{10}(4), 569--582.

\bibitem{pet}
M.-H. Qin, P. Huang, M. Zhang, C. Xiao, J.-Y. Chen, Y.-D. Zhou, M. Lu, J.-Y. Li and Y.-Q. Lan (2025). \href{https://doi.org/10.1021/jacs.5c08549}{Binuclear metal-based covalent organic framework mimicking metallohydrolases for direct photoreforming of PET plastic}. \emph{J. Am. Chem. Soc.} \textbf{147}(32), 29232--29240.
\href{https://doi.org/10.1021/jacs.5c08549.s001}{Supporting Information, Sect.~S1.5 and Table~S7}.
\bibitem{kernelpert}
B. Beckermann, M. Putinar, E. B. Saff and N. Stylianopoulos (2021).
\href{https://doi.org/10.1007/s10208-020-09458-9}{Perturbations of Christoffel--Darboux kernels: detection of outliers}.
\emph{Foundations of Computational Mathematics} \textbf{21}(1), 71--124.
\bibitem{whittall}
K.~P. Whittall and A.~L. MacKay (1989).
\href{https://doi.org/10.1016/0022-2364(89)90011-5}{Quantitative interpretation of NMR relaxation data}.
\emph{Journal of Magnetic Resonance} \textbf{84}(1), 134--152.
\bibitem{parker}
R.~L. Parker and Y.-Q. Song (2005).
\href{https://doi.org/10.1016/j.jmr.2005.03.002}{Assigning uncertainties in the inversion of NMR relaxation data}.
\emph{Journal of Magnetic Resonance} \textbf{174}(2), 314--324.
\bibitem{kaye}
B. Kaye, P.~J. Foster, T.~Y. Yoo and D.~J. Needleman (2017).
\href{https://doi.org/10.1371/journal.pone.0169337}{Developing and testing a Bayesian analysis of fluorescence lifetime measurements}.
\emph{PLoS ONE} \textbf{12}(1), e0169337.
\href{https://github.com/bryankaye1/bayesian-analysis-of-fluorescent-lifetime-data/tree/f0db816995ea2aede852547dee251d82be36492e/low_photon}{Photon histograms and instrument response}.
\bibitem{huggins}
P. Huggins, J.~S. Martin, A.~R.~J. Downey and S.~H. Won (2025).
\href{https://doi.org/10.1016/j.snb.2024.137018}{Interpretable machine learning for predicting the derived cetane number of jet fuels using compact TD-NMR}.
\emph{Sensors and Actuators B: Chemical} \textbf{426}, 137018.
\bibitem{arts}
P. Huggins, J. Martin, A. Downey and S.~H. Won. \href{https://github.com/ARTS-Laboratory/Dataset-time-domain-NMR-hydrocarbon-and-fuel-relaxometry/tree/b42b70089bf97f75e3062b1487f77cd35e85f800}{Dataset time domain NMR hydrocarbon and fuel relaxometry}. ARTS Laboratory, version 2 validation workbook, repository revision of 24 November 2026.
\bibitem{hurlimann}
W. H\"urlimann (2015).
\href{https://doi.org/10.1186/s13660-015-0709-1}{An explicit version of the Chebyshev--Markov--Stieltjes inequalities and its applications}.
\emph{Journal of Inequalities and Applications} \textbf{2015}, 192.
\bibitem{deboor}
C. de Boor (2005).
\href{https://arxiv.org/abs/math/0502036}{Divided differences}.
\emph{Surveys in Approximation Theory} \textbf{1}, 46--69.
\bibitem{denoyelle}
Q. Denoyelle, V. Duval and G. Peyr\'e (2017).
\href{https://doi.org/10.1007/s00041-016-9502-x}{Support recovery for sparse super-resolution of positive measures}.
\emph{Journal of Fourier Analysis and Applications} \textbf{23}, 1153--1194.
\bibitem{schmudgen}
K. Schm\"udgen (2017).
\href{https://doi.org/10.1007/978-3-319-64546-9}{\emph{The Moment Problem}}.
Graduate Texts in Mathematics, vol.~277, Springer, Cham.
\bibitem{stark}
P. B. Stark (1992).
\href{https://doi.org/10.1029/92JB00739}{Inference in infinite-dimensional inverse problems: discretization and duality}.
\emph{Journal of Geophysical Research} \textbf{97}, 14055--14082.
\bibitem{eftekhari}
A. Eftekhari, J. Tanner, A. Thompson, B. Toader and H. Tyagi (2021).
\href{https://doi.org/10.1016/j.acha.2019.08.004}{Sparse non-negative super-resolution---simplified and stabilised}.
\emph{Applied and Computational Harmonic Analysis} \textbf{50}, 216--280.
\bibitem{companion}
K. Zhang and M. Zhao (2026).
\href{https://doi.org/10.5281/zenodo.24756333}{Global statistical resolution of positive kinetic mixtures under correlated measurements}.
Zenodo preprint.
\bibitem{ss}
I.~J. Schoenberg and G. Szeg\H{o} (1960).
\href{https://www.numdam.org/item/CM_1959-1960__14__260_0/}{An extremum problem for polynomials}.
\emph{Compositio Mathematica} \textbf{14}, 260--268.
\end{thebibliography}
\end{document}